\documentclass{article}
\title{Orthogonal Dualities of Dynamic Stochastic Higher Spin Vertex Models, using the Drinfeld Twister} 
\author{Jeffrey Kuan and Zhengye Zhou} 
\date{}
\usepackage{amsfonts}
\usepackage{graphics}
\usepackage{float}
\graphicspath{{./xmPictures/}}

\usepackage{helvet}

\usepackage{amssymb}
\usepackage{epstopdf}

\usepackage{color}

\usepackage[utf8]{inputenc}

\usepackage{hyperref}
\usepackage{amsmath,amsthm,tikz}
\usetikzlibrary{arrows.meta}

\newtheorem{thm}{Theorem}[section]

\newtheorem{theorem}{Theorem}[section]
\newtheorem{proposition}[thm]{Proposition}

\newtheorem{corollary}[thm]{Corollary}
\newtheorem{definition}[thm]{Definition}

\theoremstyle{remark}
\newtheorem{remark}{Remark}
\newtheorem{example}{Example}

\begin{document}
{
\maketitle
}

\begin{center}\fbox{\parbox{0.8\textwidth}{ {Accessibility Statement:  A WCAG2.1AA compliant version of this PDF will be available  at \href{https://xerxes.ximera.org/jeffrey-kuan-tracywidomdrinfeldtwister/newCourse/newFolder/algebraic_part}{a Ximera webpage.}}}}

\end{center}

\begin{center}

\hspace{1in}

{Abstract}
\end{center}

{
We introduce a new algebraic method to construct duality functions for integrable dynamic models. This method will be implemented on dynamic stochastic higher spin vertex models, where we prove that the resulting duality functions between the dynamic stochastic higher spin vertex models and non-dynamic stochastic higher spin vertex models are the \( {}_3 \varphi_2\) functions. A degeneration of these duality functions is dual $q$-Krawtchouk polynomials, which  agree with the orthogonal polynomial dualities of Groenevelt–Wagenaar \cite{groenevelt2023generalized} between dynamic ASEP and  ASEP. 
}

{
The method relies on the universal twister of \(U_q(\mathfrak{sl}_2)\), regarded as a quasi-triangular quasi-Hopf algebra. Since the algebraic construction is formulated in a general setting, it is expected to produce duality functions for many other dynamic integrable models as well.

}

\normalsize


{\section{Introduction}}
{
Markov duality is a widely used tool in the theory of Markov processes, wherein an analysis of a Markov process is reduced to an analysis of a ``simpler'' Markov process. However, despite its ubiquity, discovering Markov duality  often remains something of  a ``black art''\cite{JK14}. Over the last 30 years, starting with the seminal papers \cite{SS94,Schutz97}, researchers have developed an algebraic machinery to discover Markov duality in interacting particle systems. 
}

{
In more recent years, there have been new developments on so-called ``dynamic'' models. Dynamic stochastic vertex models and dynamic ASEP were introduced by Borodin \cite{BOR17}, using the Felder-Varchenko elliptic quantum group. His result produced a contour integral formula for expectations of certain observables which resembled a duality result. Subsequent work of \cite{BC20} proved that the observable is indeed a duality function between dynamic ASEP and usual ASEP, using direct calculations involving $q$-Hermite polynomials. Later work of Aggarwal \cite{AA18} applied fusion to produce dynamic stochastic higher spin vertex models with explicit formulas for the weights, and also produced contour integral formulas for expectations. 
}

{
In this paper, the main probabilistic result is a duality for the dynamic stochastic higher spin vertex models. The duality functions are themselves vertex weights, which originated from discussions with Michael Wheeler and William Mead (private communication).
}

{
The proof utilizes algebraic machinery. In particular, the new ingredients are (1) the universal Drinfeld twister of the quantum group \(U_q( \mathfrak{sl}_2)\), in contrast to the use of the Felder-Varchenko elliptic quantum group; and (2) calculations are done using the universal $R$-matrix rather than applying fusion. As far as we know, this is the first probabilistic application of the previous two tools. The algebraic methods are expected to apply to other Drinfeld-Jimbo quantum groups.
}

We also make a remark on orthogonal polynomial duality, where the duality functions are polynomials which are orthogonal with respect to the reversible measures of the process. Because of the orthogonality, there are \textbf{a priori} applications in probability, such as to Boltzmann-Gibbs principles \cite{ACR18,ACR21} or correlation formulas \cite{CFG21,KZ23}. In this sense, orthogonal polynomial dualities are an improvement on ``triangular'' dualities, for which the probabilistic applications are not immediately obvious. A recent paper \cite{FKZ22} developed a new algebraic method for constructing orthogonal polynomial dualities in interacting particle systems, generalizing the previous methods of \cite{CFGR18,CFG21,Wolter19,Zhou_2021}. The method developed in the present paper relies on Wheeler’s observation that duality functions can themselves be interpreted as vertex weights. Extending this perspective to dualities arising from orthogonal polynomials will require new algebraic insights, which we discuss in Section~\ref{OPD}.

{
Relatedly, Groenevelt and Wagenaar constructed in \cite{groenevelt2023generalized} a dynamic ASEP$(q,j)$ admitting an orthogonal polynomial duality with respect to the reversible measures. We show that their duality function arises as a degeneration of the one obtained in the present paper.
}

{
\textbf{Acknowledgements.} The authors would like to thank (in no particular order) Michael Wheeler, William Mead, Wolter Groenevelt, Alexei Borodin, Amol Aggarwal, Ivan Corwin, Carel Wagenaar, Konstantin Matetski, and Leonid Petrov for helpful discussions. The discussions with Wolter Groenevelt and William Mead were held at the workshop ``Recent Developments in Stochastic Duality'' in Eindhoven, with support from EURANDOM and the NWO grant 613.001.753 ``Duality for interacting particle systems.''

The first author was supported by National Science Foundation grant DMS--2000331.
}

{\section{ Preliminaries and Notations}}

\subsection{Orthogonal Polynomials}\label{OPB}
In this section, we briefly review the orthogonal polynomials and special functions that will appear throughout the paper. In particular, the orthogonal dualities are related to Racah–Wilson and dual $q$-Krawtchouk polynomials, which also arise through $q-6j$ symbols.

{
Let $n$ be a nonnegative integer, the $q$-Pochhammer symbol is defined by:
\[
(a ; q)_n=\prod_{j=0}^{n-1}\left(1-q^j a\right)
\qquad
(a ; q)_{-n}=\prod_{j=1}^{n} \frac{1}{1-q^{-j} a}.
\]
From the $q$-Pochhammer symbols, we can define the  hypergeometric series
\[
{ }_p \varphi_r\left(\begin{array}{c}
a_1, a_2, \ldots, a_p \\
b_1, b_2, \ldots, b_r
\end{array} \mid q, z\right)=\sum_{k=0}^{\infty} \frac{z^k}{(q ; q)_k} \prod_{j=1}^p\left(a_j ; q\right)_k \prod_{j=1}^r\left(b_j ; q\right)_k^{-1}.
\]
For any integer $n$, define the $q$-deformed integers and factorial by
\[[n]=\frac{q^n-q^{-n}}{q-q^{-1}},\quad\quad[n]!=\frac{(q^2;q^2)_n}{(q^{-1}-q)^n}q^{-n(n+1)/2}.\]
 We extend this definition such that for any $a\in\mathbf{R}$ and  integer $n$,
\[\frac{[a+n]!}{[a]!}=\frac{(q^{2(a+1)};q^2)_n}{(q^{-1}-q)^n}q^{-na-n(n+1)/2}.\]
}

{
Now we introduce the Racah-Wilson polynomials \cite{AW79}, which is given by the $q$-hypergeometric function:
\begin{equation}
    W_n(x;\alpha,\beta,\gamma,M|q)={ }_4 \varphi_3\left(\begin{array}{c}
q^{-n}, q^{n+1}\alpha\beta, q^{-x},q^{x-M}\gamma \\
q^{}\alpha,q\beta\gamma, q^{-M} 
\end{array} q, q\right).
\end{equation}
It satisfies the following orthogonal relation:
\begin{equation}
 \sum_{x=0}^N W_n(x;\alpha,\beta,\gamma,M|q)  W_m(x;\alpha,\beta,\gamma,M|q) w(x)  =\delta_{m, n}h_n,
 \end{equation}
where

    \begin{align*}
         h_n=&\frac{(q ; q)_n(1-\alpha \beta q)(\beta q ; q)_n\left(\alpha \gamma^{-1} q ; q\right)_n\left(\alpha \beta q^{M+2} ; q\right)_n(q^{-M} \gamma )^n}{(\alpha \beta q ; q)_n\left(1-\alpha \beta q^{2 n+1}\right)(\alpha q ; q)_n(\beta \gamma q ; q)_n(q^{-M}  ; q)_n}\\
        &\times\frac{\left(q^{-M+1} \gamma  ; q\right)_{\infty}\left(\alpha^{-1} \beta^{-1} q^{-M-1} ; q\right)_{\infty}\left(\alpha^{-1} \gamma ; q\right)_{\infty}\left(\beta^{-1} ; q\right)_{\infty}}{\left(\alpha^{-1} q^{-M} \gamma;q\right)_{\infty}\left(\beta^{-1} q^{-M}  ; q\right)_{\infty}(\gamma q ; q)_{\infty}\left(\alpha^{-1} \beta^{-1} q^{-1} ; q\right)_{\infty}}, \\
        w(x)=&\frac{(q^{-M} \gamma; q)_x\left(1-q^{-M+2x} \gamma \right)(\alpha q ; q)_x(\beta \gamma q ; q)_x(q^{-M}  ; q)}{(q ; q)_x(1-q^{-M} \gamma )\left(\alpha^{-1} q^{-M} \gamma  ; q\right)_x\left(\beta^{-1} q^{-M}  ; q\right)_x(\gamma q ; q)_x(\alpha \beta q)^x}.
    \end{align*}

}

{
Next, we define the $q-6j$ symbol \cite{KR88}

    \begin{align}\label{eq: 6j}
  \notag   &\left\{\begin{array}{ccc}
a & b & e\\
d& c & f
\end{array}\right\}_q \\ \notag
 =&(-1)^{2c+2d+2e}\sqrt{[2e+1][2f+1]}\Delta\left(a,b,e\right)\Delta\left(a,c,f\right)
\Delta\left(c,e,d\right)\Delta\left(d,b,f\right)\\ \notag
& \times\sum_{w }\left(\frac{(-1)^{-w}[a+b+c+d+1-w]!}{[c+d-e-w]![b+d-f-w]![ a+c-f-w]![ a+b-e-w]!}\right.\\ 
 &\left.\times\frac{1}{[w]![w+e+f-b-c]![w-a-d+e+f]!}\right),   
    \end{align}
where the sum is taken only over  $w$ with non-negative arguments in square brackets inside the summation and
\[
\Delta\left(I_1, I_2, j_3\right)= \sqrt{\frac{\left[-I_1+I_2+j_3\right] !\left[I_1-I_2+j_3\right] !\left[I_1+I_2-j_3\right] !}{\left[I_1+I_2+j_3+1\right] !}}.
\]
The $q-6j$ symbol can be expressed through the Racah-Wilson polynomial, by setting
\[
\begin{array}{lll}
    n=a+b-e, & x=c+d-e, &M=a+b+c+d+1,\\
    \alpha=q^{2(-a-d+e+f)}, &\beta=q^{2(-a-b-c+d+1)}, &\gamma=q^{2(a+e+f-d+1)} 
\end{array}
\]
so that 
\begin{equation}\label{eq: 6j vs RW}
    \left\{\begin{array}{ccc}
a & b & e\\
d& c & f
\end{array}\right\}_q
 =\sqrt{\frac{w(x)}{h_n}}  W_n(x;\alpha,\beta,\gamma,M|q^2). 
\end{equation}
}

The dual $q$-Krawtchouk polynomial $K_n(x ; c, N \mid q)$ is defined as 
\begin{equation}
 K_n(x ; c, N \mid q)=
{ }_3 \phi_2\left(\begin{array}{c}
q^{-n}, q^{-x}, -c q^{x-N} \\
q^{-N}, 0
\end{array} ; q, q\right), 
\end{equation}
It satisfies the following orthogonal relation:
\begin{equation}
\begin{aligned}
& \sum_{x=0}^N \frac{\left(-c q^{-N}, q^{-N} ; q\right)_x}{(q, -c q ; q)_x} \frac{\left(1+c q^{2 x-N}\right)}{\left(1+c q^{-N}\right)} (-c)^{-x} q^{x(2 N-x)} K_m(x ; c, N \mid q) K_n(x ; c, N \mid q) \\
& =\left(-c^{-1} ; q\right)_N \frac{(q ; q)_n}{\left(q^{-N} ; q\right)_n}\left(c q^{-N}\right)^n \delta_{m n}, \quad c>0.
\end{aligned}
\end{equation}

\subsection{The Dynamic Stochastic Six Vertex Model}\label{model}

\begin{figure}
    \includegraphics[height=5in]{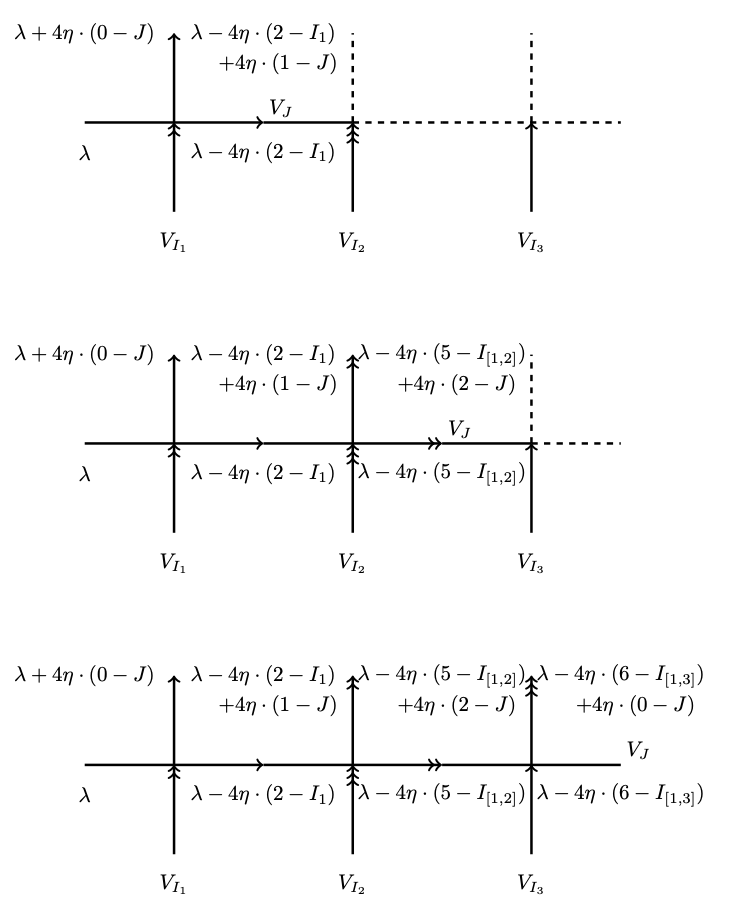}
    \centering
    \caption{This shows one-time update of the dynamic stochastic higher spin vertex model. The update occurs sequentially from left to right. The dashed lines indicate that the update has not yet occurred, and each solid line shows the number of arrows.}
    \label{fig:my_label}
\end{figure}

The purpose of this section is to introduce the dynamical stochastic
(higher spin) vertex model constructed in \cite{AA18,BOR17}.  They assign weights to configurations of up-right paths
in the square grid (See Figure \ref{fig:my_label}), and the weight of a configuration
is the product of weights of its vertices. 

 We first define the weights at a single vertex for the dynamic stochastic six vertex model:        Fix $q,\lambda,z\in \mathbb{ C}$, define
\begin{equation}\label{eq:Bweights}
\begin{gathered}
\text {weight}(\odot)=S(0,0,0,0)=1, \quad \text {weight}_{\lambda, w}(+)=S(1,1,1,1)=1, \\
\text {weight}(\mid)=S(1,0,1,0)=\frac{(q^2  -e^{-2\pi \mathbf{i} \lambda})  (z - 1)  }{ (z  q^2 - 1)  (1-e^{-2\pi \mathbf{i} \lambda})}, \\
\text {weight}(\lrcorner)=S(0,1,1,0)= \frac{(q^2 - 1)  (z-e^{-2\pi \mathbf{i} \lambda} )  }{ (z  q^2 - 1)  (1-e^{-2\pi \mathbf{i} \lambda} )}, \\
\text {weight}(\ulcorner)=S(1,0,0,1)=\frac{(q^2 - 1)  (1-e^{-2\pi \mathbf{i} \lambda}  z )  }{ (z  q^2 - 1)  (1-e^{-2\pi \mathbf{i} \lambda} )}, \\
\text {weight}(--)=S(0,1,0,1)=\frac{(z - 1)  (1-e^{-2\pi \mathbf{i} \lambda}  q^2 )  }{ (z  q^2 - 1)  (1-e^{-2\pi \mathbf{i} \lambda} )}.
\end{gathered}
\end{equation}
Here, $\lambda$ is the dynamic parameter and $z$ is the spectral parameter. 

We note that, after the change of variables \(z=q^{-1/2}w^{-1}\) and the substitution \(q\mapsto q^{1/2}\), our weights agree with those of \cite[equation~(9.7)]{BOR17}.

{
These weights were generalized by Aggarwal in \cite[Definition 1.1]{AA18} through a fusion procedure, leading to the fused weights. Let $q,\lambda,z\in \mathbb{ C}$ and $J,I\in\frac{1}{2}\mathbb{Z}_{+}$, define

\begin{equation}\label{eq: Amolweights}
\begin{aligned}
&S_{I,J,u,\lambda}\left(i_1,j_1;i_2,j_2\right) \\
&=\mathbf{1}_{i_1+j_1=i_2+j_2}
q^{2(j_2-i_1)J}
\left(\frac{u}{s}\right)^{j_1}
\frac{
(q^{2(i_1-j_2+1)};q^2)_{j_2}
(q^{2(j_2-J)};q^2)_{j_1}
}{
(su;q^2)_{i_1+j_1}
(q^2;q^2)_{j_2}
}
\\
&\times
\frac{
(suq^{2J};q^2)_{i_1-j_2}
(s^2q^{2(i_1-j_2)};q^2)_{j_1}
}{
(q^{2(j_2-J)};q^2)_{j_1-j_2}
}
\\
&\times
\frac{
(us^{-1}\kappa^{-1}q^{-2i_1};q^2)_{j_2}
(q^{2(1-i_1-J)}u^{-1}s^{-1}\kappa^{-1};q^2)_{j_1}
}{
(q^{2(1-j_2)}\kappa^{-1};q^2)_{j_1}
(q^{2(j_2-i_1-J+1)}s^{-2}\kappa^{-1};q^2)_{j_1}
}
\\
&\times
\frac{
(q^2\kappa^{-1};q^2)_{j_1}
(q^{2(j_2-2i_1+1)}s^{-2}\kappa^{-1};q^2)_{i_1-j_2}
}{
(q^{2(j_2-2i_1-J)}s^{-2}\kappa^{-1};q^2)_{j_2}
(q^{2(2j_2-2i_1-J+1)}s^{-2}\kappa^{-1};q^2)_{i_1-j_2}
}
\\
&\times
{}_{10}W_9\Bigl(
q^{-2j_2}\kappa^{-1};
q^{-2j_1},
q^{-2j_2},
q^{2(j_1-J)}\kappa^{-1},
\\
&\hspace{3.6cm}
q^{2(1-i_1)}s^{-2}\kappa^{-1},
su^{-1}q^{2(i_1-j_2+1)},
usq^{2(i_1-j_2+J)},
\\
&\hspace{3.6cm}
q^{-2i_1}\kappa^{-1};
q^2,q^2
\Bigr).
\end{aligned}
\end{equation}
 with $$u=\frac{1}{zq^2} \quad s=q^{-I} ; \quad \kappa=q^{-2(J-2 j_1)} e^{-2 \pi \mathbf{i} \lambda}$$
and where
\begin{equation}
{ }_{r+1} W_r\left(a_1 ; a_4, a_5, \ldots, a_{r+1} ; q, z\right)=\sum_{k=0}^{\infty} \frac{z^k\left(a_1 ; q\right)_k}{(q ; q)_k} \frac{1-a_1 q^{2 k}}{1-a_1} \prod_{j=4}^{r+1} \frac{\left(a_j ; q\right)_k}{\left(q a_1 / a_j ; q\right)_k}
\end{equation}
denotes the very-well poised basic hypergeometric series. 
These weights sum to $1$ in the sense that
\[
\sum_{i_{2}, j_{2}} S\left(i_{1}, j_{1} ; i_{2}, j_{2}\right)=1.
\]
In \eqref{eq: Amolweights}, $J$ is the horizontal spin and $I$ is the vertical spin. As before $\lambda$ is the dynamic parameter, and $z$ is the spectral parameter. We also note that \eqref{eq: Amolweights} agrees with that of \cite[Definition 1.1]{AA18} with substitution \(q\mapsto q^{1/2}\). 
}

{
From these weights, a discrete-time Markov chain can be constructed as a vertex model (see \cite{AA18} for a detailed definition). In this paper, we consider the model on finite sites with closed boundaries (i.e. no  horizontal incoming arrows through the left boundary and no  horizontal outgoing arrows through the right boundary), or on the infinite line.
}

{
Now we briefly describe the model defined on $[1,N]\times \mathbb{Z}_+$ (See Figure \ref{fig:my_label}).

For each vertex $(x,y) \in [1,N]\times \mathbb{Z}_+$, there is an arrow configuration, which is a quadruple \[(i_1, j_1;i_2, j_2) = (i_1, j_1;i_2, j_2)_{(x,y)} =\left(i_1(x, y), j_1(x, y);i_2(x, y), j_2(x, y)\right)\]
of non-negative integers, where $j_1$  denotes the number of horizontal incoming arrows at $(x, y)$; $i_1$ denotes the number of
vertical incoming arrows; $j_2$ denotes the number of horizontal outgoing arrows; and $i_2$ denotes the
number of vertical outgoing arrows. See Figure \ref{fig:arrow} for an example. At any vertex, the total number of incoming arrows is equal to the total number of outgoing arrows (i.e. $i_1+j_1=i_2+j_2$). Let $J_y$ be the horizontal spin  for strip $y$, and $I_x$ be the vertical spin for strip $x$, which means $0\le j_1(\cdot,y),j_2(\cdot,y) \le  2J_y$, and $0\le i_{1}(x,\cdot),i_{2}(x,\cdot)\le 2I_x$. 
\begin{figure}
\centering
 \begin{tikzpicture}[scale=1.2,>=Stealth]

\draw[fill=white] (0,0) circle (0.12);

\draw[->, thick] (-0.08,-1.2) -- (-0.08,-0.12);
\draw[->, thick] ( 0.08,-1.2) -- ( 0.08,-0.12);

\draw[->, thick] (-1.4,-0.12) -- (-0.12,-0.12);
\draw[->, thick] (-1.4, 0.00) -- (-0.12, 0.00);
\draw[->, thick] (-1.4, 0.12) -- (-0.12, 0.12);

\draw[->, thick] (0,0.12) -- (0,1.2);

\draw[->, thick] (0.12,-0.18) -- (1.4,-0.18);
\draw[->, thick] (0.12,-0.06) -- (1.4,-0.06);
\draw[->, thick] (0.12, 0.06) -- (1.4, 0.06);
\draw[->, thick] (0.12, 0.18) -- (1.4, 0.18);

\end{tikzpicture}
\caption{arrow configuration at a single vertex with $(2,3,1,4)$}
    \label{fig:arrow}
\end{figure}

In addition, each vertex $(x,y)$ is associated with a dynamical parameter $\lambda(x,y)$ according to the following rules:
\begin{align*}\lambda(x+1,y)&=\lambda(x,y)-4\eta (i_1(x+1,y)-I_{x+1}),\\ \lambda(x,y+1)&=\lambda(x,y)+4\eta (j_1(x,y)-J_{y}),\end{align*}
where $\eta$  defined by  $q=e^{2\pi \mathrm{i} \eta}$ is the asymmetry parameter. 
}

{
Given the empty left boundary (i.e. $j_1(1,\cdot)=0$ ) and arbitrary bottom boundary (i.e. $i_1(\cdot,1)$), we now update the configuration from left to right and from bottom to top with   transition probability at site $(x,y)$ defined as 
\begin{equation}
    \mathcal{P}_{x,y}(i_2,j_2|i_1,j_1)=S_{I_x,J_y,u_x,\lambda(x,y)}\left(i_1, j_1 ; i_2, j_2\right).
\end{equation}
}
{
Similarly to the stochastic six vertex model \cite{BCG16}, there is  an equivalent way to define the model as  an interacting particle system: consider the $\mathcal{P}$-distributed random
configuration and cut it by horizontal lines $y=t+\frac{1}{2}$, $t=0,1,2\ldots$, then $\eta_x(t)=i_1\left(x,t+\frac{1}{2}\right)$ is a Markov chain. Since  the horizontal arrows are directed to the right, particles may only jump to the right, and we therefore say that the model has asymmetry to the right. Equivalently, reversing the horizontal direction of the arrows yields a model with asymmetry to the left.
}

\subsection{Algebraic Background}\label{QGas}
In this section, we review the algebraic structures underlying the dynamic stochastic higher spin vertex models studied in this paper. We adopt the notation and conventions for \(U_q( \mathfrak{sl}_2)\) from \cite{BBB95}.

{
The quantum group \(U_q( \mathfrak{sl}_2)\) is a unital algebra generated by \(E_{\pm}\) and \(q^H\) with relations:
\[
\left[H, E_{ \pm}\right]= \pm 2 E_{ \pm}, \quad\left[E_{+}, E_{-}\right]=\frac{q^{H}-q^{-H}}{q-q^{-1}}.
\]
The co-unit is $\epsilon(q^H)=1,\epsilon(E_{\pm})=0$. The co-product in this paper will be
\[
\Delta(H)=H \otimes \mathrm{id}+\mathrm{id} \otimes H, \quad \Delta\left(E_{ \pm}\right)=E_{ \pm} \otimes q^{H / 2}+q^{-H / 2} \otimes E_{ \pm}.
\]
The above definitions define a bialgebra structure on \(U_q( \mathfrak{sl}_2)\). There is also an antipode which defines a Hopf algebra structure, but we will not (explicitly) use the antipode.
}

{In the present paper, the universal R-matrix will ultimately give rise to the stochastic transition weights of the  models. We now recall Drinfeld's universal $R$-matrix:
\[
R_{12}^{D}=q^{\frac{1}{2} H \otimes H} \sum_{i=0}^{\infty}\left(q-q^{-1}\right)^{i} \frac{q^{-\frac{i(i+1)}{2}}}{[i] !} q^{\frac{i}{2} H} E_{+}^{i} \otimes q^{-\frac{i}{2} H} E_{-}^{i}.
\]
This $R$-matrix satisfies the property that
\[
R^D \Delta( X ) = \Delta'(X) R^{D}
\]
for any \(X\in U_q( \mathfrak{sl}_2)\), where \(\Delta'=P \circ \Delta\) is the opposite comultiplication with $P$ defined by \(P(X\otimes Y) = Y\otimes X\). This property makes \(U_q( \mathfrak{sl}_2)\) an \textit{almost-cocommutative} Hopf algebra. The $R$-matrix additionally satisfies
\[
R^{D}_{12}R^{D}_{21}=1, \quad (\epsilon \otimes \epsilon)R^{D}_{12}=1 \quad R^{D}_{12}(\Delta \otimes \mathrm{id})(R^{D}) = R^{D}_{23}(\mathrm{id} \otimes \Delta)(R^{D}),
\]
which makes it a \textit{coboundary} Hopf algebra. Additionally, 
\[
(\Delta \otimes \mathrm{id})(R^{D}) = R^{D}_{13}R^{D}_{23}, \quad (\mathrm{id} \otimes \Delta)(R^{D}) = R^{D}_{13}R^{D}_{12},
\]
making it a \textit{quasi-triangular} Hopf algebra. In any quasi-triangular Hopf algebra, the $R$-matrix satisfies the quantum Yang-Baxter equation \cite{KSBOOK}
\[
R^{D}_{12}R^{D}_{13}R^{D}_{23} = R^{D}_{23}R^{D}_{13}R^{D}_{12}.
\]
}

{
In fact, there is a spectral-dependent version of the $R$-matrix, denoted \(R^D(z)\),  satisfying the  Yang-Baxter Equation
\begin{equation}
R^{D}_{12}(u)R^{D}_{13}(uv)R^{D}_{23}(v) = R^{D}_{23}(v) R^{D}_{13}(uv) R^{D}_{12}(u).
\end{equation}
}

{
There is a generalization of Hopf algebras to quasi-Hopf algebras, which can be constructed with a twisting procedure. Suppose there is a twister satisfying the following shifted cocycle condition \cite{BBB95}:
\begin{equation}\label{eq: cocycle}
F_{12}(\lambda)(\Delta\otimes \mathrm{id} )F(\lambda)=F_{23}({\lambda-2\eta H_1})(\mathrm{id}  \otimes\Delta)F(\lambda);
\end{equation}
then the  twisted $R$-matrix 
\begin{equation}\label{eq: twistR}
  R(z,\lambda) = F_{21}(\lambda) \left(R^D(z)\right)^{-1} F_{12}^{-1}(\lambda)  
\end{equation}
satisfies the dynamical Yang-Baxter Equation \cite{FV96}

\begin{equation}\label{YBE}
 \begin{aligned}  & R_{12}(u, { \lambda - 2\eta H_3})R_{13}(uv,{ \lambda})R_{23}(v,{\lambda -2\eta H_1}) \\
&\quad\quad\quad\quad\quad\quad= R_{23}(v,{ \lambda}) R_{13}(uv,{ \lambda -2\eta H_2}) R_{12}(u, { \lambda}),
\end{aligned}
\end{equation}
 where $\eta$ is defined in a way such that \(q=e^{2\pi \mathbf{i} \eta}\). 
}

{
Now set 
\[\Delta_{\lambda}(X)=F(\lambda)\Delta(X)F^{-1}(\lambda),\]
which is called the twisted co-product. Then,
\[
R(z,\lambda) \Delta_{\lambda}( X ) = \Delta_{\lambda}'(X) R(z,\lambda)
\]
for any \(X\in U_q( \mathfrak{sl}_2)\). Additionally, 
\[
(\Delta_{\lambda} \otimes \mathrm{id})(R(z,\lambda)) =\Phi_{312} R_{13}(z,\lambda)\Phi^{-1}_{132} R_{23}(z,\lambda)\Phi_{123},
\]
\[ (\mathrm{id} \otimes \Delta_{\lambda})(R(z,\lambda)) = \Phi^{-1}_{231} R_{13}(z,\lambda)\Phi_{213}R_{12}(z,\lambda)\Phi^{-1}_{123}
\]
where
\(\Phi=\left(F_{23}(\lambda)(\mathrm{id} \otimes \Delta) F(\lambda)\right) \left(F_{12}(\lambda)(\Delta \otimes \mathrm{id}) F(\lambda)\right)^{-1}\). These properties make \(U_q( \mathfrak{sl}_2)\) a quasi-Hopf  algebra. For more details about quasi-Hopf  algebras, we refer the readers  to \cite{JKOS97}.
}

{

There is also a $*$-Hopf algebra structure on \(U_q( \mathfrak{sl}_2)\) for nonzero real values of $q$. The involution is given by the explicit formula
\begin{equation}\label{eq: *structure}
E_{+}^{*}=E_{-} , \quad E_{-}^{*}= E_{+}, \quad\left(q^{H}\right)^{*}=q^{H}.
\end{equation}

{
Finally, recall that for \(J\in \{1/2,1,3/2,2,5/2,\ldots\}\), the $(2J+1)$-dimensional representations are given by
\begin{align*}
\rho^{(J)}(H)|J+ m\rangle & =2 m|J+ m\rangle &  -J\le m\le J \\
\rho^{(J)}\left(E_{ \pm}\right)|J+ m\rangle & =\sqrt{[J \mp m][J \pm m+1]}|J+ m \pm 1\rangle,
\end{align*}
Denote this representation by $V^{2J}$. These representations provide the local state spaces for the higher spin vertex models studied later. Under this interpretation, we interpret $J+m$ to be the number of particles (arrows) and $J-m$ to be the number of holes at each site. 
}



{\section{Main Results}}

{
The main theorems of this paper consist of several Markov duality results, which can be stated independently of the algebraic framework. }

Let us start with the definition of Markov duality
\begin{definition}
Two Markov processes \(\xi(t)\) and $\eta(t)$ on states spaces $\mathfrak{X}$ and $\mathfrak{Y}$, respectively, are \textit{dual} with respect to the function $D(\xi,\eta)$ on $\mathfrak{X} \times \mathfrak{Y}$ if
\[
\mathbf{E}_\xi[D(\xi(t),\eta)] = \mathbf{E}_\eta[D(\xi,\eta(t))]\quad \text{ for all } t\geq 0.
\]
If, in addition, $\mathfrak{X}=\mathfrak{Y}$ and $\xi(t)=\eta(t)$ then $\xi(t)$ is \textit{self-dual} with respect to the function $D(\xi,\eta)$.
\end{definition}
}

{
If $\mathfrak{X}$ and $\mathfrak{Y}$ are discrete, an equivalent definition of duality is the intertwining relation
$$
L_\xi D = DL_\eta^T,
$$
where $L_\xi$ is the generator of $\xi(t)$ viewed as a matrix with rows and columns indexed by $\mathfrak{X}$, $L_\eta$ is the generator of $\eta(t)$ with rows and columns indexed by $\mathfrak{Y}$, and $D$ is the duality function viewed as a matrix with rows indexed by $\mathfrak{X}$ and columns indexed by $\mathfrak{Y}$. The superscript $ ^T$ indicates matrix transpose. Note that this is using the convention in probability, that a stochastic matrix has rows summing to $1$.
}

{
Recall that $I_i$ denotes the vertical spin at site $i$. Let $\mu$ and $\xi$ denote the Markov processes associated with the vertex models whose transition matrices are given by $\mathcal{T}^{\mathrm{stoch}}(z,\lambda)$ and $\mathcal{T}^{\mathrm{stoch}}_{\mathrm{rev}}(z)$, defined in \eqref{T1} and \eqref{T2}, respectively. Then the following duality results hold as a consequence of Theorem~\ref{Multi Site Duality}.
}

\begin{theorem}\label{thm: main}
The dynamic stochastic higher spin vertex model $\mu$, with asymmetry to the right, is dual to a (non-dynamic) stochastic higher spin vertex model $\xi$ with asymmetry to the left, with respect to the duality functions:
    \begin{align*}
   \mathcal{D}_c(\mu,\xi)=\prod_{i=1}^N&\left(e^{-2\pi \mathbf{i}\lambda}q^{2(-c+N_{[1,i-1]}(\mu+\xi-2I)+\mu_i-2I_i)};q^2\right)_{2I_i-\mu_i}\\
  &\times { }_3 \varphi_2\left(\begin{array}{c}
 q^{-2(2I_i-\mu_i)},q^{-2\xi_i} ,e^{-2\pi \mathbf{i}\lambda}q^{2(2N_{[1,i-1]}(\mu-I)+\mu_i-2I_i)}\\
q^{-4I_i},e^{-2\pi \mathbf{i}\lambda}q^{2(-c+N_{[1,i-1]}(\mu+\xi-2I)+\mu_i-2I_i)}\end{array} ; q^2, q^{2}\right) \\
&\times  q^{  -4I_iN_{[1,i-1]}(\xi)-2\xi_iN_{[1,i]}(\mu)},
    \end{align*}
where $c$ is an arbitrary  number in $\mathbb{C}$, $I$ is the vector $(I_1,\cdots,I_N)$ and $N_{[1,i]}(\eta)=\sum_{k=1}^i\eta_k$.

\end{theorem}

{
Fix $q$, and take a limit of $c$ such that  $q^c \to 0.$  In this limit, the duality becomes orthogonal:

\begin{corollary}\label{cor: ort}
    The dynamic stochastic higher spin vertex model $\mu$, with asymmetry to the right, is dual to a (non-dynamic) stochastic higher spin vertex model $\xi$ with asymmetry to the left, with respect to the duality functions:
 \begin{align*}\label{Dort}
     \mathcal{D}_{ort}(\mu,\xi)= \prod_{i=1}^N & { }_3 \varphi_2\left(\begin{array}{c}
 q^{2(2I_i-\mu_i)},q^{2\xi_i} ,e^{-2\pi \mathbf{i}\lambda}q^{2(2N_{[1,i-1]}(\mu-I)+\mu_i)}\\
q^{4I_i},0\end{array} q^{-2}, q^{-2}\right) \\ &\times q^{ -4I_iN_{[1,i]}(\xi)-2\xi_iN_{[1,i-1]}(\mu)}
\end{align*}

Moreover, there are orthogonal relations
\begin{align}
& \sum_{\xi} \mathcal{D}_{ort}(\mu,\xi)\mathcal{D}_{ort}(\mu',\xi) w(\xi ; q, I)=\frac{\delta_{\mu, \mu^{\prime}}}{W(\mu ; q, I ,\lambda)}, \\
& \sum_{\mu} \mathcal{D}_{ort}(\mu,\xi)\mathcal{D}_{ort}(\mu,\xi')W(\mu ; q, I ,\lambda)=\frac{\delta_{\xi, \xi^{\prime}}}{w(\xi ; q, I)} ,
\end{align}
with 
\[
    w(\xi ; q, I)=q^{\frac{1}{2}\sum_{i=1}^N\left(2I_i\xi_i-2\xi_iN_{[1,i]}(2I)\right)}\prod_{i=1}^Nq^{-\xi_i(\xi_i-2I_i)}\left[\begin{array}{c}
2I_i \\
\xi_i
\end{array}\right]_{q^{-2}}
\]
and
\[
 W(\mu ; q, I ,\lambda)= \prod_{i=1}^N \mathcal{W}\left(2I_i-\mu_i ; q^{-1}, 2I_i, \log_{q^{-2}}\left(-e^{-2\pi \mathbf{i}\lambda}q^{4N_{[1,i-1]}(\mu-I)}\right)\right),
\]

\begin{equation}
\mathcal{W}(x ; q, N, \rho) =\frac{1+q^{4 x+2 \rho-2 N}}{1+q^{2 \rho-2 N}} \frac{\left(-q^{2 \rho-2 N} ; q^2\right)_x}{\left(-q^{2 \rho+2} ; q^2\right)_x} \frac{q^{-x(2 \rho+1+x-2 N)}}{\left(-q^{-2 \rho} ; q^2\right)_N}\left[\begin{array}{c}
N \\
x
\end{array}\right]_{q^2} .
\end{equation}
\end{corollary}
}

{
\begin{remark}
    $ \mathcal{D}_{ort}(\mu,\xi)$ is equal to $K_L$ in   \cite[Theorem 3.17]{groenevelt2023generalized}, where  $K_L$ is proven to be an orthogonal duality between dynamic ASEP and ASEP, the orthogonal relation is given by  \cite[Corollary 3.18]{groenevelt2023generalized}. The measure $w$ and $W$ are reversible measures for ASEP and dynamic ASEP, respectively. If we let $I_i=I$ for any $i$, then
    up to a constant, \begin{equation*}
        \lim_{I\xrightarrow[]{}\infty}\mathcal{D}_{ort}(\mu,\xi)=\prod_{i=1}^N q^{2\xi_iN_{[1,i]}(\mu)},
    \end{equation*}which is a self duality function for $q$-TAZRP (totally asymmetric zero range process).
\end{remark}
}

{
In the limit to the non-dynamic model as $\lambda=-\mathbf{i}\infty$ in  $ \mathcal{T}^{stoch}(z,\lambda)$ (i.e. $e^{-2\pi\mathbf{i}\lambda} \to \infty$ ), we obtain a new duality function:
\begin{corollary}
    The  stochastic higher spin vertex model, with asymmetry to the right, is dual to a stochastic higher spin vertex model with asymmetry to the left, with respect to the duality function
    \begin{align*}
        \mathcal{D}_{\text{new}}(\mu,\xi)
 =\prod_{i=1}^N&\frac{\left(q^{2(C_0+N_{[1,i-1]}(\mu+\xi-2I)+\mu_i-2I_i+1)};q^2\right)_{2I_i-\mu_i}}{\left(q^{-2(C_0+N_{[1,i-1]}(\mu+\xi-2I)+\xi_i)};q^2\right)_{2I_i-\mu_i}}\\
  &\times { }_2 \varphi_1\left(\begin{array}{c}
 q^{-2(2I_i-\mu_i)},q^{-2(2I_i-\xi_i)} \\
q^{-4I_i}\end{array} q^2,q^{-2(C_0+N_{[1,i]}(\mu+\xi-2I))} \right)\\
&\times q^{  
  -4I_iN_{[1,i-1]}(\xi)-2\xi_iN_{[1,i]}(\mu)},
    \end{align*}
where \(C_0\) is an arbitrary fixed constant.

\end{corollary}
}

{

\begin{remark}

    Take $C_0\xrightarrow[]{}\infty$ in $\mathcal{D}_{\text{new}}$ we get a triangular duality
\begin{equation*}
    \begin{aligned}
        \mathcal{D}_{tr}(\mu,\xi)=&\prod_{i=1}^N  1_{\{\mu_i\ge \xi_i\}}\frac{[\mu_i]![2I_i-\xi_i]!}{\left[\mu_i-\xi_i\right]!}q^{-4I_iN_{[1,i-1]}(\xi)  +\mu_i\xi_i+2\mu_iN_{[1,i-1]}(\xi)-2I_i\xi_i}.    
    \end{aligned}
\end{equation*}
This triangular duality coincides with the duality in \cite[Theorem 2.5(a)]{Kuan18}. Indeed, applying the ``charge-parity symmetry'' to $\xi$ in Theorem 2.5(a) gives the same result.

\end{remark}
}

{
A simple way to check the duality result is to look at ``two-site'' dualities or degenerate cases. In the next examples, we display a few ``two-site'' duality matrices which are calculated up to a constant in the total number of particles.  Below, $ \mathcal{T}^{stoch}(z,\lambda)$ is the ``two-site''  transition matrix for process $\eta$ and $\mathcal{T}^{\mathrm{stoch}}_{\text{rev}}(z)$ is the transition matrix for the non-dynamic process $\xi$.
\[
\mathcal{T}^{\mathrm{stoch}}_{\text{rev}}(z)= \begin{bmatrix}
 1&                     0&                         0& 0\\
0& \frac{q^2 - 1}{(z  q^2 - 1)}& \frac{q^2  (z - 1)}{z  q^2 - 1}& 0\\
0&   \frac{z - 1}{z  q^2 - 1}& \frac{z  (q^2 - 1)}{z  q^2 - 1}& 0\\
0&                     0&                         0& 1\\
 
    \end{bmatrix},
\]

\[
  \mathcal{T}^{stoch}(z,\lambda)= \begin{bmatrix}
  1&                                           0&                                           0& 0  \\
  
   0&  \frac{(q^2 - 1)  (1-e^{-2\pi \mathbf{i} \lambda}  z )  }{ (z  q^2 - 1)  (1-e^{-2\pi \mathbf{i} \lambda} )}&  \frac{(q^2  -e^{-2\pi \mathbf{i} \lambda})  (z - 1)  }{ (z  q^2 - 1)  (1-e^{-2\pi \mathbf{i} \lambda})} & 0  \\ 0&   \frac{(z - 1)  (1-e^{-2\pi \mathbf{i} \lambda}  q^2 )  }{ (z  q^2 - 1)  (1-e^{-2\pi \mathbf{i} \lambda} )}&    \frac{(q^2 - 1)  (z-e^{-2\pi \mathbf{i} \lambda} )  }{ (z  q^2 - 1)  (1-e^{-2\pi \mathbf{i} \lambda} )}& 0  \\
   0&                                           0&                                           0& 1  \\
    \end{bmatrix}.    
\]

}

{
\begin{example} 
The matrix \(\mathcal{D}_c\) equals

\resizebox{1\linewidth}{!}{
$
\begin{bmatrix}
(e^{-2\pi\mathbf{i}\lambda}q^{-2c-2}-1)(e^{-2\pi\mathbf{i}\lambda}q^{-2c-4}-1)&-(e^{-2\pi\mathbf{i}\lambda}(e^{-2\pi\mathbf{i}\lambda}-q^{2c+2})(q^{2c}-1))q^{-4c-6}&-(e^{-2\pi\mathbf{i}\lambda}(e^{-2\pi\mathbf{i}\lambda}-q^{2c+2})(q^{2c}-1))q^{-4c-6}&(e^{-4\pi\mathbf{i}\lambda}(q^{2c}-q^2)(q^{2c}-1))q^{-4c-8}\\
1-e^{-2\pi\mathbf{i}\lambda}q^{-2c-2}&-(e^{-2\pi\mathbf{i}\lambda}q^{-2c-2}-1)q^{-2}&(e^{-2\pi\mathbf{i}\lambda}(q^{2c}-1))q^{-2c-4})&(e^{-2\pi\mathbf{i}\lambda}(q^{2c}-1))q^{-2c-6}\\
1-e^{-2\pi\mathbf{i}\lambda}q^{-2c-2}&(e^{-2\pi\mathbf{i}\lambda}(q^{2c+2}-1))q^{-2c-4}&-(q^{-2c}e^{-2\pi\mathbf{i}\lambda}-1)q^{-4}&(e^{-2\pi\mathbf{i}\lambda}(q^{2c}-1))q^{-2c-6}\\
1&q^{-4}&q^{-4}&q^{-8}\\
    \end{bmatrix}
$}
and satisfies
$ \mathcal{T}^{\mathrm{stoch}}(z, { \lambda}) \mathcal{D}_c =   \mathcal{D}_c \left(\mathcal{T}^{\mathrm{stoch}}_{\text{rev}}(z)\right)^*$.
\end{example}
}

{
\begin{example} The matrix $\mathcal{D}_{\text{new}}$ equals
\small
\[
\mathcal{D}_{\text{new}}= \begin{bmatrix}
((q^{2C_0}-q^2)(q^{2C_0}-q^4))q^{-6}&(q^{2C_0}(q^{2C_0}-q^2))q^{-6}&(q^{2C_0}(q^{2C_0}-q^2))q^{-6}&q^{4C_0-6}\\
1-q^{2C_0-2}&-(q^{2C_0}-q^2)q^{-4}&-q^{2C_0-4}&-q^{2C_0-6}\\
1-q^{2C_0-2}&-q^{2C_0-4}&-(q^{2C_0}-1)q^{-4}&-q^{2C_0-6}\\
1&q^{-4}&q^{-4}&q^{-8}\\
    \end{bmatrix}
\]
and satisfies
\(\mathcal{T}^{\mathrm{stoch}}(z, {-\mathbf{i}\infty }) \mathcal{D}_{\text{new}} =   \mathcal{D}_{\text{new}} \left(\mathcal{T}^{\mathrm{stoch}}_{\text{rev}}(z)\right)^*\).
\end{example}
}

{
\begin{example}
The matrix $\mathcal{D}_{\text{tr}}$ equals
\[
\mathcal{D}_{tr}= \begin{bmatrix}
       1&         0&       0&     0\\
q^{-1/2}&       1&       0&     0\\
  q^{-1/2}&         0&       q^{-2}&     0\\
        q& q^{3/2}& q^{3/2}& q^4\\ 
    \end{bmatrix}
\]
and satisfies \(\mathcal{T}^{\mathrm{stoch}}(z, {-\mathbf{i}\infty }) \mathcal{D}_{tr} =   \mathcal{D}_{tr} \left(\mathcal{T}^{\mathrm{stoch}}_{\text{rev}}(z)\right)^*\).
\end{example}
}

{\section{Algebraic Framework}\label{OPD}}
The proof of the duality relations is based on a general algebraic framework. We begin by describing this framework and then mention that the dynamic stochastic higher-spin vertex model satisfies its assumptions.

{
Consider the following setup: Suppose   the matrix \(R(z,\lambda)\) acts on a tensor product \(V\otimes W\), which depends on a quantization parameter $q$, a spectral parameter $z$ and the dynamical parameter \(\lambda\) and satisfies the  dynamical Yang-Baxter Equation \eqref{YBE}. We also assume that \(R(z,\lambda)\) is obtained from twisting universial $R$ matrix with a twister $F(\lambda)$ satisfying \eqref{eq: cocycle}:  \begin{equation}\label{eq: ass1}
    R(z,\lambda) = F_{21}(\lambda) R(z) F_{12}^{-1}(\lambda).
\end{equation}   

We consider the tensor product \(V_0 \otimes V_1 \otimes \cdots \otimes V_N\), with \(V_0\) called the ``auxiliary space'' and \(V_1 \otimes \cdots \otimes V_N\) called ``the bulk space''. Let $2J_0+1$ be the dimension of $V_0$ and  \(2I_i+1\) be the dimension of \(V_i\) ($i>0$). We write $R_{i,j}$ to denote the action of $R$ on the tensor product $V_i\otimes V_j$.

{
Let $P$ be the permutation map on two sites, defined by \(P: V\otimes W \rightarrow W \otimes V\) with \(P(v\otimes w)=w\otimes v\). Let the \(\check{ \ }\) and \ \(\hat{ }\) notation refer to 
\[
\check{R}_{i,i+1} = P_{i,i+1}R_{i,i+1} = R_{i+1,i}P_{i,i+1}.
\]
\[
\hat{R}_{i,i+1} = R_{i,i+1}P_{i,i+1} =P_{i,i+1} R_{i+1,i}.
\]
 Define 
\begin{multline*}
\mathcal{B}(z, {\lambda}) = R_{0,N}\left(\infty, \lambda - 2\eta N_{[1,N-1]}(H)\right) R_{0,N-1}\left(\infty,\lambda - 2\eta N_{[1,N-2]}(H) \right) \\
\cdots R_{0,1}\left(\infty,\lambda\right),
\end{multline*}
where \(N_{[a,b]}(H)=H_a+H_{a+1}\cdots+H_{b}\). For each \(0\le m,l \leq \dim V_0\), we denote
\[
\langle \mu | B_{ml}(z,\lambda) | \xi \rangle = \langle m, \mu | \mathcal{B}(z,\lambda) | l,\xi \rangle.
\]
}

{
\begin{theorem}\label{Intertwining Theorem}
Suppose that the matrix $R(z)$ in \eqref{eq: ass1}  satisfies the following transposition formula:
\begin{equation}\label{Transposition}
 \Pi^{\otimes 2} R(z) = R(z)^* \Pi^{\otimes 2}
\end{equation}
with \(\Pi\) being a matrix such that 
\begin{equation} \label{eq: PIP}
  \Pi^{\otimes 2} P =P \Pi^{\otimes 2}.  
\end{equation}

Then there is the intertwining relationship
\begin{equation}\label{eq: interwining}
  \check{R}_{i,i+1}(z,\lambda-2\eta N_{[1,i-1]}(H))   B_{ml}(z,\lambda)  = B_{ml}(z,\lambda){\mathfrak{R}}_{i,i+1}(z,\lambda(l))^*,
\end{equation}
where $1\le i \le N-1$ and
\[
    {\mathfrak{R}}(z,\lambda)^*:=F_{12}(\lambda) (\Pi^{-1})^{\otimes 2} (F_{21}(\lambda))^*\hat{R}(z,\lambda)^*
  (F_{21}(\lambda)^*)^{-1}\Pi^{\otimes 2}  F_{12}^{-1}(\lambda),
\]
\(\lambda(l)= \lambda -2\eta N_{[1,i-1]}(H)- 2\eta H_0(l)\) with \(H_0(l)\) denoting the constant \(2(l-J_0).\)
\end{theorem}
}

{
If a few more properties hold, then there is an intertwining relationship between  stochastic matrices. First suppose there is a ground state symmetry that turns $R(z,\lambda)$ into a stochastic matrix $S(z,\lambda)$:
\begin{equation}\label{GroundState}
G^{-1}_{i,i+1}(\lambda)\check{R}_{i,i+1}(z,\lambda)G_{i,i+1}(\lambda) = \check{S}_{i,i+1}(z,\lambda)
\end{equation}
such that \(S_{i,i+1}(z,\lambda)\) is a stochastic matrix, meaning that each row sums to $1$ and all entries are non-negative  for some values of \(z,\lambda\).

Last, suppose there exists a stochastic matrix $\mathfrak{S}(z)$ and another ground state matrix $\Tilde{G}$ such that
\begin{equation}\label{eq: combined}
  \hat{R}(z,-\mathbf{i}\infty)^*=\Tilde{G} \hat{\mathfrak{S}}(z)^*  \Tilde{G}^{-1}. 
\end{equation}

{
Notice that 
   \begin{align*}
    \lim_{H_0(l)\xrightarrow{}\infty}{\mathfrak{R}}(z,\lambda(l))^*&= (\Pi^{-1})^{\otimes 2}  \hat{R}(z,-\mathbf{i}\infty)^*\Pi^{\otimes 2}, 
\end{align*}
which means it reduces to the non-dynamic case in this limit. In the theorem below, we will take the quantity $H_0(l)$ to $\infty$, as in such limit, the twister terms in ${\mathfrak{R}}(z,\lambda)$ become identities.
\begin{theorem}\label{Two Site Duality}
With all of the above assumptions (specifically equations \eqref{YBE}, \eqref{eq: ass1}-\eqref{eq: PIP}, \eqref{GroundState}-\eqref{eq: combined}, there is the ``two-site'' duality relation 

\begin{equation*}
    \check{S}_{i,i+1}\left(z,\lambda-2\eta N_{[1,i-1]}(H)\right)D_{i,i+1}=D_{i,i+1} \hat{\mathfrak{S}}(z)^*_{i,i+1},
\end{equation*}
where 
\begin{equation*}
    D_{i,i+1}=\text{const}\cdot G_{i,i+1}^{-1}(\lambda-2\eta N_{[1,i-1]}(H))\left(\lim_{H_0(l)\xrightarrow{}\infty}B_{ml}(z,\lambda) \right)\left(\Pi^{-1}\right)^{\otimes 2} \Tilde{G}_{i,i+1}
\end{equation*}

\end{theorem}
}
\begin{remark}
    On the right-hand-side of the duality relation, the  term ``const'' is a normalization constant so that the duality function is nontrivial in the limit.
\end{remark}

\begin{remark}
The assumptions \eqref{eq: ass1}-\eqref{eq: PIP} are expected to apply quite broadly to dynamical $R$-matrices arising from Drinfeld twists of quantum groups. In particular, these assumptions are natural in the quasi-Hopf algebra setting and are not specific to $U_q(\mathfrak{sl}_2)$. 
Furthermore, in many highest-weight representations, the vacuum vector is an eigenvector of the $R$-matrix; this naturally leads to the ground-state transformations  \eqref{GroundState}-\eqref{eq: combined}.  We therefore expect the framework to extend to other stochastic models arising from $R$-matrices of quasi-triangular quasi-Hopf algebra.
\end{remark}

Now we extend this duality relation to finitely many sites. Define the stochastic transfer matrix acting on the bulk space by 
\begin{equation}\label{T1}
   \begin{aligned}
  \mathcal{T}^{\mathrm{stoch}}(z,{\lambda}) &=\check{S}_{1,2}(z,\lambda)\check{S}_{2,3}(z,\lambda-2\eta H_1)\cdots \\
 &\cdots \check{S}_{N-2,N-1}\left(z,\lambda - 2\eta N_{[1,N-3]}(H) \right)\check{S}_{N-1,N}\left(z, \lambda - 2\eta N_{[1,N-2]}(H)\right) 
\end{aligned} 
\end{equation}
and
\begin{equation}
    \label{T2}\mathcal{T}^{\mathrm{stoch}}_{\text{rev}}(z) 
=\hat{\mathfrak{S}}_{N-1,N}\hat{\mathfrak{S}}_{N-2,N-1}\cdots  \hat{\mathfrak{S}}_{1,2}  . 
\end{equation}  
}
We also define the corresponding non-stochastic transfer matrix by
\begin{equation}\label{T0}
   \begin{aligned}
  \mathcal{T}(z,{\lambda}) =&\check{R}_{1,2}(z,\lambda)\check{R}_{2,3}(z,\lambda-2\eta H_1)\cdots \\
 & \cdots \check{R}_{N-2,N-1}\left(z,\lambda - 2\eta N_{[1,N-3]}(H) \right)\check{R}_{N-1,N}\left(z, \lambda - 2\eta N_{[1,N-2]}(H)\right) \end{aligned} 
\end{equation}
and
\begin{equation}
    \label{T3}\mathcal{T}^{}_{\text{rev}}(z) 
=\hat{R}_{N-1,N}(z,-\mathbf{i}\infty)\hat{R}_{N-2,N-1}(z,-\mathbf{i}\infty)\cdots  \hat{R}_{1,2}(z,-\mathbf{i}\infty)  . 
\end{equation}  
{
\begin{remark}
   Notice that, in the definition of $\mathcal{T}^{\mathrm{stoch}}(z,\lambda)$, each $\check{S}$ is stochastic, i.e., its rows sum to $1$. Therefore, $\mathcal{T}^{\mathrm{stoch}}(z,\lambda)$ may be viewed as the transfer matrix of the dynamical stochastic vertex model with asymmetry to the right, where the dynamical parameter at site $1$ is equal to $\lambda$. See Figure~\ref{fig:my_label} for the corresponding updates of the parameters. Similarly, $\mathcal{T}^{\mathrm{stoch}}_{\text{rev}}$ can be viewed as the transfer matrix of the non-dynamic stochastic vertex model with asymmetry to the left. The subscript in $\mathcal{T}^{\mathrm{stoch}}_{\text{rev}}$ refers to space reversal.
\end{remark}
}

{
\begin{theorem}\label{Multi Site Duality}
With the same assumptions as in Theorem \ref{Intertwining Theorem}, we further assume that there exist ground-state transformations $\mathcal{G}(\lambda)$ and $\Tilde{\mathcal{G}}$ such that
\begin{equation}\label{eq: multiGS}
    \mathcal{T}^{\mathrm{stoch}}(z,{\lambda})=\mathcal{G}^{-1}(\lambda) \mathcal{T}(z,{\lambda})\mathcal{G}(\lambda),
\end{equation}
and 
\begin{equation}\label{eq: multiGSrev}
     \left(\mathcal{T}^{}_{\text{rev}}(z)\right)^*=\Tilde{\mathcal{G}} \left(\mathcal{T}^{\mathrm{stoch}}_{\text{rev}}(z)\right)^*  \Tilde{\mathcal{G}}^{-1}.
\end{equation}
Then there is the duality relation 
\[
\mathcal{T}^{\mathrm{stoch}}(z, { \lambda}) \mathcal{D} =   \mathcal{D} \left(\mathcal{T}^{\mathrm{stoch}}_{\text{rev}}(z)\right)^* ,
\]
where 
\[
\mathcal{D}= \text{const}\cdot
\mathcal{G}^{-1}(\lambda)\left( \lim_{H_0(l)\xrightarrow{}\infty} B_{ml}(z,\lambda)\right)\left(\Pi^{-1}\right)^{\otimes N} \Tilde{\mathcal{G}}.
\]

\end{theorem}
}
\begin{remark}
    On the right-hand-side of the duality relation, the  term ``const'' is a normalization constant so that the duality function is nontrivial in the limit.
\end{remark}

To relate these algebraic constructions to the dynamic stochastic higher-spin vertex model, we will verify that this model fits within the general framework developed above. The following propositions establish this fact.

\begin{proposition}\label{prop: equations satisfied}
The dynamic stochastic higher-spin vertex model satisfies assumptions made in Theorem \ref{Two Site Duality} and \ref{Multi Site Duality}.
\end{proposition}

The next step is to identify explicitly the duality functional $\mathcal{D}$ from Theorem~\ref{Multi Site Duality} and show that it coincides with the expression stated in Theorem~\ref{thm: main}.
\begin{proposition}\label{prop: calc}
  The duality functional $\mathcal{D}$ has the expression stated in Theorem \ref{thm: main}.
\end{proposition}
The final step is to show that the fusion procedure of Aggarwal \cite{AA18} yields the same stochastic model as the construction based on the universal $R$-matrix and twister developed in this paper. The following proposition establishes this identification.

\begin{proposition}\label{prop: matching}
The stochastic weights \eqref{eq: Amolweights} from \cite{AA18} coincide with the weights $S(z,\lambda)$.
\end{proposition}

{\section{Proofs of Algebraic Framework}}

\subsection{Proof of Theorem \ref{Intertwining Theorem}}
This section is devoted to the proof of Theorem \ref{Intertwining Theorem}. We first derive an intertwining relation from the dynamical Yang--Baxter equation, then translate it into a statement for the matrix elements $B_{ml}(z,\lambda)$, and finally use the transposition formula \eqref{Transposition} to complete the proof.

{\subsubsection{Yang-Baxter-Equation}}

{Replacing the lattice sites $\{1,2,3\}$ with $\{0,i,i+1\}$ in the dynamical Yang-Baxter equation \eqref{YBE}, we get
\begin{multline*}
R_{0,i}(u,  \lambda - 2\eta H_{i+1}) R_{0,i+1}(uv,  \lambda) R_{i,i+1}(v,  \lambda - 2\eta H_0) \\
= R_{i,i+1}(v, {\lambda}) R_{0,i+1}(uv, {\lambda-2\eta H_i}) R_{0,i}(u, { \lambda}).
\end{multline*}
Now letting $P_{i,j}$ be the permutation operator that permutes site $i$ and $j$, and since $P^2$ is the identity,
\begin{align*}
&P_{i,i+1} R_{0,i}(u,{ \lambda - 2\eta H_{i+1}}) P_{i,i+1}P_{i,i+1} R_{0,i+1}(uv,{ \lambda }) P_{i,i+1}P_{i,i+1} R_{i,i+1}(v,{ \lambda - 2\eta H_{0}}) P_{i,i+1}\\
&= P_{i,i+1} R_{i,i+1}(v, { \lambda}) R_{0,i+1}(uv, { \lambda-2\eta H_{i}}) R_{0,i}(u, { \lambda}) P_{i,i+1}.
\end{align*}
Now using that
$$
P_{j,k} R_{i,j}(z, { \lambda}   ) P_{j,k} = R_{i,k}(z, { \lambda}), \quad { P_{j,k} R_{ij}(z, { \lambda-2\eta H_{k}}) P_{j,k} = R_{i,k}(z, { \lambda-2\eta H_{j}})  },
$$
we get 
\begin{multline*}
  R_{0,i+1}(u, { \lambda-2\eta H_i}) R_{0,i}(uv, { \lambda}) R_{i+1,i}(v, {\lambda-2\eta H_0})\\
= P_{i,i+1} R_{i,i+1}(v, { \lambda}) R_{0,i+1}(uv, { \lambda-2\eta H_i}) R_{0,i}(u, { \lambda}) P_{i,i+1}.  
\end{multline*}
Plugging in
$$
\check{R}_{i,i+1}(z, {\lambda} ) := P_{i,i+1}R_{i,i+1}(z, { \lambda} ) = R_{i+1,i}(z, { \lambda})P_{i,i+1},
$$
we obtain
\begin{multline*}
R_{0,i+1}(u, { \lambda- 2\eta H_{i}}) R_{0,i}(uv, { \lambda}) \check{R}_{i,i+1}(v, { \lambda-2\eta H_0})\\
=  \check{R}_{i,i+1}(v, { \lambda}) R_{0,i+1}(uv, { \lambda-2\eta H_i}) R_{0,i}(u, { \lambda}) .
\end{multline*}
Notice that if we take the limit as $u\rightarrow\infty$ (or let $u=0$, we don't consider this case in the current paper) then 
$$\lim_{u\xrightarrow{}{\infty}}R_{0,i+1}(u, { \lambda- 2\eta H_{i}}) R_{0,i}(uv, { \lambda}) =\lim_{u\xrightarrow{}{\infty}} R_{0,i+1}(uv, { \lambda-2\eta H_i}) R_{0,i}(u, { \lambda}) .$$

Thus, for any $i=1,\ldots,N-1$,
$$
\mathcal{B}(z,\lambda) \check{R}_{i,i+1}(q,\lambda-2\eta N_{[1,i-1]}(H)-2\eta H_0) = \check{R}_{i,i+1}(q,\lambda-2\eta N_{[1,i-1]}(H)) \mathcal{B}(z,\lambda).
$$
}

\subsubsection{Reducing from $\mathcal{B}$ to $B$}

{
We continue from the last section, recall that for each $m,l \leq \dim V_0$
$$
\langle \mu | B_{ml}(z,\lambda) | \xi \rangle = \langle m, \mu | \mathcal{B}(z,\lambda) | l,\xi \rangle,
$$
we have that
\begin{multline*}
    \left\langle m,\mu \vert   \check{R}_{i,i+1}(z,\lambda-2\eta N_{[1,i-1]}(H))  \mathcal{B}(z,\lambda) \vert l,\xi  \right\rangle =\\ \left\langle m,\mu \vert \mathcal{B}(z,\lambda) \left(  \check{R}_{i,i+1}(z, \lambda-2\eta N_{[1,i-1]}(H) - 2\eta H_0) \right)  \vert l,\xi  \right\rangle
\end{multline*}
with $m+N(\mu)=l+N(\xi)$, where $N(\mu)$ is the total number of particles in $\mu$.
This implies
\begin{multline*}
   \sum_{\nu} \langle \mu \vert   \check{R}_{i,i+1}(z,\lambda-2\eta N_{[1,i-1]}(H)) \vert \nu \rangle \langle m, \nu | \mathcal{B}(z,\lambda) | l,\xi \rangle  \\ = \sum_{\xi} \langle m, \mu | \mathcal{B}(z,\lambda) | l, \xi\rangle \langle l, \xi |  \check{R}_{i,i+1}(z, \lambda -2\eta N_{[1,i-1]}(H)-2\eta H_0)  | l,\xi \rangle, 
\end{multline*}
which can be stated as
\begin{multline*}
 \check{R}_{i,i+1}(z,\lambda-2\eta N_{[1,i-1]}(H))   B_{ml}(z,\lambda) \\
 = B_{ml}(z,\lambda) \check{R}_{i,i+1}(z,\lambda-2\eta N_{[1,i-1]}(H)-2\eta H_0(l))
\end{multline*}
for all values of $m$ and $l$. Here, $H_0(l)$ is the constant $2(l-J_0)$.
}

{\subsubsection{Transposition}}
{

Now, by \eqref{eq: ass1} we get
\begin{equation*}
 (F_{12}(\lambda)^*)^{-1}\Pi^{\otimes 2} F_{21}(\lambda)^{-1} {R}(z,\lambda) = R(z,\lambda)^* (F_{21}(\lambda)^*)^{-1} \Pi^{\otimes 2}F_{12}(\lambda)^{-1}.
\end{equation*}
Next, multiplying $P$ on both sides from left and also write $F_{21}$ on the left hand side by $PF_{21}P$ yields
\begin{equation*}
 P(F_{12}(\lambda)^*)^{-1}\Pi^{\otimes 2} PF_{12}(\lambda)^{-1}P {R}(z,\lambda) = PR(z,\lambda)^* (F_{21}(\lambda)^*)^{-1} \Pi^{\otimes 2}F_{12}(\lambda)^{-1}.
\end{equation*}
Using fact that $P$ intertwines with $\Pi^{\otimes 2}$ gives 
\begin{equation*}
 (F_{21}(\lambda)^*)^{-1}\Pi^{\otimes 2} F_{12}(\lambda)^{-1} \check{R}(z,\lambda) = \hat{R}(z,\lambda)^* (F_{21}(\lambda)^*)^{-1} \Pi^{\otimes 2}F_{12}(\lambda)^{-1}.
\end{equation*}
Notice that by definition,
\begin{equation*}
    {\mathfrak{R}}^*(z,\lambda)=\check{R}(z,\lambda).
\end{equation*}
Therefore we can write
\begin{equation*}
    \check{R}_{i,i+1}(z,\lambda-2\eta N_{[1,i-1]}(H))   B_{ml}(z,\lambda) 
 = B_{ml}(z,\lambda) {\mathfrak{R}}^*(z,\lambda(l)),
\end{equation*}
where $\lambda(l)= \lambda -2\eta N_{[1,i-1]}(H)- 2\eta H_0(l)$. 
}

\subsection{Proof of Theorem \ref{Two Site Duality} and \ref{Multi Site Duality}}
The proof of Theorem \ref{Two Site Duality} and \ref{Multi Site Duality} starts from combining \eqref{eq: interwining} with \eqref{GroundState} and then take the limit as $H_0(l)\xrightarrow[]{}\infty$, which results in 

\begin{equation}
\begin{aligned}
&  \check{S}_{i,i+1}(z,\lambda-2\eta N_{[1,i-1]}(H))  G_{i,i+1}^{-1} (\lambda-2\eta N_{[1,i-1]}(H))\left(\lim_{H_0(l)\xrightarrow{}\infty}B_{ml}(z,\lambda) \right) \\
  &= G_{i,i+1}^{-1} (\lambda-2\eta N_{[1,i-1]}(H))\left(\lim_{H_0(l)\xrightarrow{}\infty}B_{ml}(z,\lambda) \right)(\Pi^{-1})^{\otimes 2}\hat{R}(z,-\mathbf{i}\infty )^*\Pi^{\otimes 2}.   
\end{aligned}
\end{equation}
Now we plug in \eqref{eq: combined}, 
\begin{equation}
\begin{aligned}
&  \check{S}_{i,i+1}(z,\lambda-2\eta N_{[1,i-1]}(H))  G_{i,i+1}^{-1} (\lambda-2\eta N_{[1,i-1]}(H))\left(\lim_{H_0(l)\xrightarrow{}\infty}B_{ml}(z,\lambda) \right) \\
  &= G_{i,i+1}^{-1} (\lambda-2\eta N_{[1,i-1]}(H))\left(\lim_{H_0(l)\xrightarrow{}\infty}B_{ml}(z,\lambda) \right)\\
  &\quad\quad\times(\Pi^{-1})^{\otimes 2}\Tilde{G} {\hat{\mathfrak{S}}}(z)^*  \Tilde{G}^{-1}\Pi^{\otimes 2}.\\
\end{aligned}
\end{equation}

To finish the proof of Theorem \ref{Multi Site Duality}, we need to extend from a single time update to multiple time updates. These are standard arguments, which we repeat below. First, to interpret the intertwining of transfer matrices as a particle system satisfying duality, we apply the same arguments as in section 4.4 of \cite{Kuan17}. Also, note that the dynamic parameter  changes in each time step. Since updating each time step is equivalent to fusing the horizontal spin  and the dynamic parameter changes in fusion (section 4.3 of \cite{AA18}, our duality result extends  to multiple time steps.

{\section{Proofs of Application to \(U_q( \mathfrak{sl}_2)\)}}

{\subsection{Proof of Proposition \ref{prop: equations satisfied}}}
{
In this section, we show that the universal $R$-matrix of \(U_q( \mathfrak{sl}_2)\) introduced in section \ref{QGas} satisfies the transposition formula \eqref{Transposition} and that there exist ground-state symmetries \eqref{GroundState} and \eqref{eq: multiGS}. We also construct a stochastic matrix $\mathfrak{S}(z)$ together with its  ground-state matrices $\Tilde{G}$ and $\Tilde{\mathcal{G}}$   such that \eqref{eq: combined}  and \eqref{eq: multiGSrev} hold.

}

\subsubsection{Transposition Formula}
{
First of all, there is a transposition formula for $R^D(z)$ \cite[equation (12)]{KMMO16} 
$$\Pi^{\otimes 2}\Lambda^{\otimes 2}R^D(z)= R^D(z)^*\Pi^{\otimes 2}\Lambda^{\otimes 2},$$
where $\Pi$ is a permutation matrix. More specifically,
\begin{equation}\label{eq: Pi}
  \langle \beta|\Pi|\alpha\rangle=1_{\{\beta=2J-\alpha\}}.  
\end{equation}
We mention here that $\Pi^2=\mathrm{id}$ and \eqref{eq: PIP} is satisfied.
Since we used a different representation from \cite{KMMO16}, we calculate $\Lambda$ here.
}
\begin{proposition}
   {$\Lambda=\mathrm{id} $ and $\Pi^{\otimes 2}R^D(z)= R^D(z)^*\Pi^{\otimes 2}.$}
\end{proposition}
\textbf{Proof.}{
As shown in \cite{KMMO16}, $\Lambda$ does not depend on $z$, thus we take $z=0$ and explicitly solve for $\Lambda$. 
}

{
Recall
$$
R_{12}^{D}(0)=q^{\frac{1}{2} H \otimes H} \sum_{i=0}^{\infty}\left(q-q^{-1}\right)^{i} \frac{q^{-\frac{i(i+1)}{2}}}{[i] !} q^{\frac{i}{2} H} E_{+}^{i} \otimes q^{-\frac{i}{2} H} E_{-}^{i}.
$$
Then applying the $*$ operator \eqref{eq: *structure}, we obtain
$$
(R_{12}^{D}(0))^*= \sum_{i=0}^{\infty}\left(q-q^{-1}\right)^{i} \frac{q^{-\frac{i(i+1)}{2}}}{[i] !} E_{-}^iq^{\frac{i}{2} H}  \otimes E_+^iq^{-\frac{i}{2} H}q^{\frac{1}{2} H \otimes H}. 
$$
}

{
When $i=\eta_1-\xi_1=\xi_2-\eta_2\ge 0$
\begin{equation*}
    \begin{aligned}
         \frac{\langle \eta_1,\eta_2|\Lambda^{\otimes 2}|\eta_1,\eta_2\rangle}{\langle \xi_1,\xi_2|\Lambda^{\otimes 2}|\xi_1,\xi_2\rangle}
    &=\frac{\langle\eta_1,\eta_2|R^D(0)|\xi_1,\xi_2\rangle}{\langle  2I_1-\eta_1,2I_2-\eta_2|R^D(0)^*|2I_1-\xi_1,2I_2-\xi_2\rangle}\\
    &=q^{2(\eta_1-I_1)(\eta_2-I_2)-2(I_1-\xi_1)(I_2-\xi_2)+i(\eta_1-\eta_2-\xi_2+\xi_1-2I_1+2I_2)}
    \\&=q^{(\eta_1-2I_1+\xi_1)(\eta_1+\eta_2-\xi_1-\xi_2)}=1.
    \end{aligned}
\end{equation*}
Thus $\langle \beta|\Lambda|\alpha\rangle=1_{\{\alpha=\beta\}} $.}

{\subsubsection{The Ground State Transformations}\label{GST}}

{
From the form of the universal $R$-matrix in section \ref{QGas}, it is immediate that
\begin{equation}\label{GS}
 \check{R}_{i,j}(z,\lambda)\vert 0,0 \rangle  = \vert 0,0 \rangle.
\end{equation}
Taking $N=2$ in $B_{ml}$, we obtain
\begin{equation*}
    \begin{aligned}
 &\check{R}_{i,i+1}(z,\lambda)  B_{ml}(z,\lambda)\vert 0,0 \rangle =  B_{ml}(z,\lambda)\check{R}_{i,i+1}(z,\lambda-2\eta h^0(l))\vert 0,0 \rangle\\
 &= B_{ml}(z,\lambda) \vert 0,0 \rangle.       
    \end{aligned}
\end{equation*}
In words, $B_{ml}(z,\lambda)  \vert 0,0 \rangle$ is an eigenvector of $ \check{R}_{i,i+1}(z,\lambda)$ of eigenvalue $1$.
}

{
Consider the value
\begin{equation}\label{eq: GSformula}
  \langle \xi \vert G_{i,i+1}(\lambda) \vert \xi\rangle =\langle \xi| B_{ml}(z,\lambda) \vert 0,0 \rangle,  
\end{equation}
where for weight conservation we require $m+\xi_1+\xi_2=l$. 
Now note that
\begin{align*}
& \sum_\mu \langle\xi_{i+1},\xi_i|\check{S}_{i,i+1}(z,\lambda)|\mu_i,\mu_{i+1}\rangle \\
& = \sum_\mu \frac{\langle\mu_i,\mu_{i+1}| G_{i,i+1}|\mu_i,\mu_{i+1}\rangle}{\langle \xi_i,\xi_{i+1}| G_{i,i+1}|\xi_i,\xi_{i+1}\rangle}\langle\xi_{i+1},\xi_i|\check{R}_{i,i+1}(z,\lambda)|\mu_i,\mu_{i+1}\rangle\\
&=\frac{ \langle\xi_{i+1},\xi_i|\check{R}_{i,i+1}(z,\lambda)   B_{ml}(z,\lambda)\vert 0,0 \rangle}{ \langle \xi_{i+1},\xi_i| B_{ml}(u,z\lambda) \vert 0,0 \rangle}\\
&=1,
\end{align*}
demonstrating that the ground state transformation exists.
}
We show the explicit formula for $G$ here, and prove it in section \ref{sec: GS}
\begin{proposition}\label{prop: G}
 Suppose   $\dim V_i=2I_i+1$ where $i\ge 1$. Then $G:V_1\otimes V_2\xrightarrow[]{}V_1\otimes V_2$ is given by
 \begin{equation}
     \begin{aligned}
        & \langle \xi_1,\xi_{2} \vert G(\lambda) \vert \xi_1,\xi_{2} \rangle\\
    & =\frac{\left([2I_1-\xi_1]![2I_2-\xi_{2}]![\xi_{2}]![\xi_1]!\right)^{-1/2}}{\left(q^{4(j(\lambda)+1+\xi_1-I_1)};q^2\right)_{2I_1-\xi_1}\left(q^{4(j(\lambda)+1+\xi_1+\xi_{2}-I_1-I_2)};q^2\right)_{2I_2-\xi_{2}}}q^{(I_1+2I_2)\xi_1+I_2\xi_2}.
     \end{aligned}
 \end{equation}
 \end{proposition}

We now extend this construction to multiple sites by defining $\mathcal{G}$ so that \eqref{eq: multiGS} holds. The proof is postponed to Section~\ref{sec: multiG}.

\begin{proposition}\label{prop: multiG}
{
Define $$ \langle \mu|\mathcal{G}(\lambda)|\mu\rangle =\prod_{i=1}^N\frac{q^{  
  2I_iN_{[1,i-1]}(\mu)+I_i\mu_i}}{\sqrt{[\mu_i]![2I_i-\mu_i]!}\left(e^{-2\pi \mathbf{i}\lambda}q^{4N_{[1,i]}(\mu-I)+2};q^2\right)_{2I_i-\mu_i}}.$$
Then
   \begin{equation*}
    \mathcal{T}^{\mathrm{stoch}}(z,{\lambda})=\mathcal{G}^{-1}(\lambda) \mathcal{T}(z,{\lambda})\mathcal{G}(\lambda).
\end{equation*}
}
\end{proposition}

{\subsubsection{Stochastic Matrix $\mathfrak{S}(z)$ }}
In this section we construct the stochastic transition matrix $\mathfrak{S}(z)$ for the dual process. we first show that the following inversion formula holds for the stochastic matrix $S(z,\lambda)$:

\begin{proposition}
\begin{equation}\label{Inversion}
\Pi^{\otimes 2}S(z,\lambda) \Pi^{\otimes 2} = \mathfrak{S}( z,\lambda),
\end{equation}
where \(\mathfrak{S}( z,\lambda)\) is defined by substituting $q$ with \(q^{-1}\) and $z$ with $z^{-1}$ in \(S(z,\lambda)\) and $\Pi$ is defined in \eqref{eq: Pi}.
\end{proposition}

\textbf{Proof.}
{
We verify the inversion formula \eqref{Inversion} directly for the spin $\frac{1}{2}$ case:
\begin{align*}
S(1,0,1,0)_{q,z\xrightarrow[]{}q^{-1},z^{-1}}=S(0,1,0,1)
\end{align*}
and also
\begin{align*}
S(0,1,1,0)_{q,z\xrightarrow[]{}q^{-1},z^{-1}}=S(1,0,0,1).
\end{align*}
Thus, this proves the inversion formula in the spin $1/2$ case. For the general case, the standard fusion argument applies (e.g. section 4.3 of \cite{AA18} or section 3.2 of \cite{Kuan17}) Q.E.D.
}

Now combine \eqref{GroundState} and \eqref{Inversion} and use the facts that $G(\lambda)$ is diagonal, $\Pi^*=\Pi$, $P^*=P$, $\Pi^{\otimes 2}P=P\Pi^{\otimes 2}$ and  $\Pi^2=P^2=\mathrm{id}$,
\begin{equation}
\begin{aligned}
 {R}_{i,i+1}(z,\lambda) &= P_{i,i+1}G_{i,i+1}(\lambda) P_{i,i+1}{S}_{i,i+1}(z,\lambda) G^{-1}_{i,i+1}(\lambda) \\
   &= P_{i,i+1}G_{i,i+1}(\lambda) P_{i,i+1}\Pi^{\otimes 2}\mathfrak{S}_{i,i+1}( z,\lambda) \Pi^{\otimes 2} G^{-1}_{i,i+1}(\lambda). 
\end{aligned} 
\end{equation}
Thus
\begin{equation}
\begin{aligned}
   \hat{R}_{i,i+1}(z,\lambda)^*&= P_{i,i+1}{R}_{i,i+1}(z,\lambda)^*\\
    &=  P_{i,i+1} G^{-1}_{i,i+1}(\lambda) \Pi^{\otimes 2}\mathfrak{S}_{i,i+1}( z,\lambda)^*\Pi^{\otimes 2} P_{i,i+1}  G_{i,i+1}(\lambda)P_{i,i+1}\\
    &=P_{i,i+1} G^{-1}_{i,i+1}(\lambda) \Pi^{\otimes 2}P_{i,i+1}\hat{\mathfrak{S}}_{i,i+1}( z,\lambda)^*\Pi^{\otimes 2} P_{i,i+1}  G_{i,i+1}(\lambda)P_{i,i+1},
\end{aligned}   
\end{equation}
we see that \eqref{eq: combined} is satisfied with $\Tilde{G}:=P_{i,i+1} G^{-1}_{i,i+1}(\lambda) \Pi^{\otimes 2}P_{i,i+1}.$

To further simplify $\Tilde{G}$, we use the following result, which will be proven  in section \ref{sec: GS}

\begin{proposition}\label{prop: PG}
    There is a diagonal matrix $C$ such that 
\begin{equation}\label{eq: PG}
  \Pi^{\otimes 2} G(-\mathbf{i}\infty)\Pi^{\otimes 2}=C P G(-\mathbf{i}\infty) P.  
\end{equation}
\end{proposition}
Using \eqref{eq: PG}, we can simply let  $\Tilde{G}:=\Pi^{\otimes 2} G^{-1}(-\mathbf{i}\infty)$. Extending this construction to multiple sites, we find that \eqref{eq: multiGSrev} holds with $\Tilde{\mathcal{G}}=\Pi^{\otimes N}\mathcal{G}^{-1}(-\mathbf{i}\infty)$.

\subsection{Proof of Proposition \ref{prop: calc}}

{
Let $\sigma'_i=\mu_i-I_i$ and $\sigma_i=\xi_i-I_i$ ($-I_i\le \sigma_i,\sigma_i'\le I_i$). By an automorphism relating the conventions of the present paper to those of \cite{BBB95}, we see that the action of $R(\infty,\lambda)$ on $V^{2I_1}\otimes V^{2I_2}$ with representations defined in section \ref{QGas} is given by the following matrix entries:

\begin{equation}\label{eq: Rformula}
\begin{aligned}
&\left\langle I_1+ \sigma_1^{\prime}\left|\left\langle I_2+\sigma_2^{\prime}\left| R^{2I_1,2I_2}(\infty,\lambda)\right| I_1+ \sigma_1\right\rangle\right| I_2+\sigma_2\right\rangle\\
&=(-1)^{\sigma_2^{\prime}-\sigma_2} q^{C_{j(\lambda)}+C_{j(\lambda)+\sigma_1+\sigma_2}-C_{ j(\lambda)+\sigma_1}-C_{j(\lambda)+\sigma_2^{\prime}}} \\
&\times\frac{\mathcal{N}_\psi^{\left(I_2\right)}\left(q^{2 j(\lambda)+1}, \sigma_2^{\prime}\right) \mathcal{N}_\psi^{\left(I_1\right)}\left(q^{2 j(\lambda)+1} q^{2 \sigma_2^{\prime}}, \sigma_1^{\prime}\right)}{\mathcal{N}_\psi^{\left(I_2\right)}\left(q^{2 j(\lambda)+1} q^{2 \sigma_1}, \sigma_2\right) \mathcal{N}_\psi^{\left(I_1\right)}\left(q^{2 j(\lambda)+1}, \sigma_1\right)}\left\{\begin{array}{ccc}
I_1 & j(\lambda)+\sigma_1+\sigma_2 & j(\lambda)+\sigma_2^{\prime} \\
I_2 & j(\lambda) & j(\lambda)+\sigma_1
\end{array}\right\}_q,
\end{aligned}
\end{equation}
where  $q^{2 j(\lambda)+1}:=e^{-\pi \mathbf{i}\lambda}$,  $C_n=n(n+1)$, the last symbol represents the  6-j coefficient \eqref{eq: 6j} (also see eq. (2.54) in \cite{CGR94}) and
\begin{multline*}
\mathcal{N}_\psi^{\left(I_1\right)}\left(q^{2 j(\lambda)+1}, \sigma_1\right)\\
=(-1)^{-2I_1+\frac{\sigma_1}{2}} \frac{\sqrt{\left[I_1+2j(\lambda)+\sigma_1+1\right] !}}{(q^{4j(\lambda)+4};q^2)_{I_1+\sigma_1}}\frac{\left(q-q^{-1}\right)^{I_1} q^{(2 j(\lambda)+1)I_1} q^{I_1 \sigma_1}}{ \sqrt{\left[2 j(\lambda)+1+2\sigma_1\right]\left[-I_1+2j(\lambda)+\sigma_1\right] !}}.
\end{multline*}
}

{
For simplicity, denote
$$f_\mathcal{N}:=\frac{\mathcal{N}_\psi^{\left(I_2\right)}\left(q^{2 j(\lambda)+1}, \sigma_2^{\prime}\right) \mathcal{N}_\psi^{\left(I_1\right)}\left(q^{2 j(\lambda)+1} q^{2 \sigma_2^{\prime}}, \sigma_1^{\prime}\right)}{\mathcal{N}_\psi^{\left(I_2\right)}\left(q^{2 j(\lambda)+1} q^{2 \sigma_1}, \sigma_2\right) \mathcal{N}_\psi^{\left(I_1\right)}\left(q^{2 j(\lambda)+1}, \sigma_1\right)}$$
and 
$$f_{6j}:=(-1)^{\sigma_2^{\prime}-\sigma_2} q^{C_{j(\lambda)}+C_{j(\lambda)+\sigma_1+\sigma_2}-C_{ j(\lambda)+\sigma_1}-C_{j(\lambda)+\sigma_2^{\prime}}}\left\{\begin{array}{ccc}
I_1 & j(\lambda)+\sigma_1+\sigma_2 & j(\lambda)+\sigma_2^{\prime} \\
I_2 & j(\lambda) & j(\lambda)+\sigma_1
\end{array}\right\}_q.$$
}

\subsubsection{Calculating Limits of $B_{ml}(z,\lambda)$}

{
In this section, we calculate several limits of $B_{ml}(z,\lambda)$ that will be used throughout the paper.
}

\begin{proposition}\label{prop: limits}
    Suppose $\dim V_0=2J_0+1$ and  $\dim V_i=2I_i+1$ when $i\ge 1$, $c=2J_0-l$ is a fixed positive integer. Let $I=(I_1,\ldots,I_N)$, then up to a constant in $N(\xi)$ and $N(\mu)$, we have the following limits:
 {   \begin{enumerate}
     \item \begin{equation*}
\begin{aligned}
   \lim_{J_0\xrightarrow[]{}\infty}&\langle \mu | B_{2J_0-c+N(\xi-\mu),2J_0-c}(\infty,\lambda) | \xi \rangle\\
 =\prod_{i=1}^N&\frac{\left(q^{2(2j(\lambda)-c+N_{[1,i-1]}(\mu+\xi-2I)+\mu_i-2I_i+1)};q^2\right)_{2I_i-\mu_i}}{\left(q^{4\left(j(\lambda)+1+N_{[1,i]}(\mu-I)\right)};q^2\right)_{2I_i-\mu_i}\sqrt{[\mu_i]![2I_i-\xi_i]![2I_i-\mu_i]![\xi_i]!}}\\
  &\times{ }_3 \varphi_2\left(\begin{array}{c}
 q^{-2(2I_i-\mu_i)},q^{-2\xi_i} ,q^{2(2j(\lambda)+2N_{[1,i-1]}(\mu-I)+\mu_i-2I_i+1)}\\
q^{-4I_i},q^{2(2j(\lambda)-c+N_{[1,i-1]}(\mu+\xi-2I)+\mu_i-2I_i+1)}\end{array} q^2, q^{2}\right)\\
& \times q^{  
  2I_iN_{[1,i-1]}(\mu-\xi)+I_i(\mu_i+\xi_i)-2\xi_iN_{[1,i]}(\mu)},
\end{aligned}
\end{equation*}
     \item \begin{equation*}
\begin{aligned}
    \lim_{J_0\xrightarrow[]{}\infty\atop j(\lambda)\xrightarrow[]{}\infty}&\langle \mu | B_{2J_0-c+N(\xi-\mu),2J_0-c}(\infty,\lambda) | \xi \rangle\\
=\prod_{i=1}^N&1_{\{\mu_i\ge \xi_i\}}\frac{1}{\left[\mu_i-\xi_i\right]!}\sqrt{\frac{[\mu_i]![2I_i-\xi_i]!}{[2I_i-\mu_i]![\xi_i]!}} q^{2I_iN_{[1,i-1]}(\mu-\xi)  +I_i(\mu_i-\xi_i)+\mu_i\xi_i+2\mu_iN_{[1,i-1]}(\xi)}, 
  \end{aligned}
\end{equation*}
     \item \begin{equation*}
\begin{aligned}
  \lim_{j(\lambda)\xrightarrow[]{}\infty\atop
  2j(\lambda)-c\xrightarrow[]{}C_0}& \lim_{J_0\xrightarrow[]{}\infty }\langle \mu | B_{2J_0-c+N(\xi-\mu),2J_0-c}(\infty,\lambda) | \xi \rangle\\
    =\prod_{i=1}^N&\frac{\left(q^{2(C_0+N_{[1,i-1]}(\mu+\xi-2I)+\mu_i-2I_i+1)};q^2\right)_{2I_i-\mu_i}}{\left(q^{4\left(j(\lambda)+1+N_{[1,i]}(\mu-I)\right)};q^2\right)_{2I_i-\mu_i}\sqrt{[\mu_i]![2I_i-\xi_i]![2I_i-\mu_i]![\xi_i]!}}\\
  &\times { }_3 \varphi_2\left(\begin{array}{c}
 q^{-2(2I_i-\mu_i)},q^{-2\xi_i} ,0\\
q^{-4I_i},q^{2(C_0+N_{[1,i-1]}(\mu+\xi-2I)+\mu_i-2I_i+1)}\end{array} q^2, q^{2}\right)\\
& \times q^{  
  2I_iN_{[1,i-1]}(\mu-\xi)+I_i(\mu_i+\xi_i)-2\xi_iN_{[1,i]}(\mu)},\\ 
  \end{aligned}
\end{equation*}
     \item \begin{equation*}
\begin{aligned}
 \lim_{c\xrightarrow[]{}\infty }  \lim_{J_0\xrightarrow[]{}\infty}&\langle \mu | B_{2J_0-c+N(\xi-\mu),2J_0-c}(\infty,\lambda) | \xi \rangle\\
&=\prod_{i=1}^N\frac{[2I_i]! q^{2I_iN_{[1,i-1]}(\mu+\xi)  +I_i(\mu_i+\xi_i)} }{\sqrt{[\mu_i]![2I_i-\xi_i]![2I_i-\mu_i]![\xi_i]!}\left(q^{4\left(j(\lambda)+1+N_{[1,i]}(\mu-I)\right)};q^2\right)_{2I_i-\mu_i}}.\\
    \end{aligned}
\end{equation*}
 \end{enumerate}}
\end{proposition}

\textbf{Proof.}
{We only show the proof of the first limit and note that the remaining limits can be derived by the same argument. 

{
We will repeatedly make use of the following fact in the proof:  as $N\rightarrow\infty$
$$
\frac{\left[\alpha \pm N\right] !}{\left[\beta \pm N\right] !} \sim(\mp)^{\alpha-\beta} \frac{q^{\mp \frac{1}{2}(\alpha-\beta)(\alpha+\beta+1)}}{\left(q-q^{-1}\right)^{\alpha-\beta}} q^{-(\alpha-\beta) N}.
$$

}

We consider the limit of \eqref{eq: Rformula} as $I_1\xrightarrow[]{}\infty$ with $I_1-\sigma_1=c$, where $c$ is a fixed integer such that $-I_i\le \sigma_i,\sigma_i'\le I_i$. Denote $c'=I_1-\sigma_1'=c+\sigma_2'-\sigma_2$. Then
\begin{equation}
    \begin{aligned}
     \lim_{I_1\xrightarrow[]{}\infty \atop I_1-\sigma_1=c}f_\mathcal{N} =& \frac{(-1)^{1/2\sigma_2+1/2\sigma_2'-I_2}}{(q-q^{-1})^{1/2\sigma_2+1/2\sigma_2'-I_2}\left(q^{4j(\lambda)+4+4\sigma_2'};q^2\right)_{I_2-\sigma_2'}} q^{I_2 - \sigma_2/4 - (3  \sigma_2')/4 + I_2  \sigma_2'}\\
   &\times q^{  2  I_2  j(\lambda)  - j(\lambda)  \sigma_2 - j(\lambda)  \sigma_2' + (c  \sigma_2)/2 + (c  \sigma_2')/2 - (\sigma_2  \sigma_2')/2 - \sigma_2^2/4 - \sigma_2'^2/4}\\  &\times\frac{\sqrt{[-c+2j(\lambda)]![\sigma_2'+I_2+2j(\lambda)+1]!}}{\sqrt{[2j(\lambda)+1+2\sigma_2'][2\sigma_2'+2j(\lambda)-c']![-I_2+\sigma_2'+2j(\lambda)]!}},     
    \end{aligned}
\end{equation}
and
{\small
\begin{equation}
    \begin{aligned}
        \lim_{I_1\to\infty \atop I_1-\sigma_1=c} f_{6j}
=&
(-1)^{4j(\lambda)+\sigma_2'-\sigma_2}
q^{\sigma_2- I_2^2/2 - I_2/2 + I_2\sigma_2 - \sigma_2'/2 + 2I_1\sigma_2} \\
&\times q^{c(-I_2-2\sigma_2) + j(\lambda)(\sigma_2-\sigma_2')+\sigma_2^2-(\sigma_2')^2/2}
\sqrt{[2j(\lambda)+2\sigma_2'+1]}\,\Delta(I_2,j(\lambda),j(\lambda)+\sigma_2') \\
&\times \sqrt{[I_2+\sigma_2]![I_2-\sigma_2]![2j(\lambda)-c]![c]![c']![2j(\lambda)+2\sigma_2'-c']!} \\
&\times \sum_{w}
\left(
\frac{(-1)^{-w}q^{w(I_2+2j(\lambda)+\sigma_2'+1)}}
{[c-w]![I_2-\sigma_2'-w]![I_2+\sigma_2-w]![w]![w+\sigma_2'-\sigma_2]![w+2j(\lambda)+\sigma_2'-I_2-c]!}
\right).
    \end{aligned}
\end{equation}

}
}

{
Thus
\begin{multline*}
    \lim_{I_1\xrightarrow[]{}\infty \atop I_1-\sigma_1=c}\left\langle I_1+ \sigma_1^{\prime}\left|\left\langle I_2+ \sigma_2^{\prime}\left| R^{2I_1,2I_2}(\infty,\lambda)\right| I_1+ \sigma_1\right\rangle\right| I_2+ \sigma_2\right\rangle\\
    =C_1(I_1,I_2,\sigma_2,\sigma_2') \mathcal{F}(I_2,j(\lambda),c,\sigma_2,\sigma_2'),
\end{multline*}
where
\normalsize
\begin{equation*}
    \begin{aligned}
        &\mathcal{F}(I_2,j(\lambda),c,\sigma_2,\sigma_2')= \frac{[-c+2j(\lambda)]!}{\left(q^{4j(\lambda)+4+4\sigma_2'};q^2\right)_{I_2-\sigma_2'}[2j(\lambda)+\sigma_2'-I_2-c]![\sigma_2'-\sigma_2]!}\\
  &\times\sqrt{\frac{[I_2+\sigma_2']![I_2-\sigma_2]![c']!}{[I_2-\sigma_2']![I_2+\sigma_2]![c]!}} { }_3 \varphi_2\left(\begin{array}{c}
q^{-2c}, q^{-2(I_2-\sigma_2')},q^{-2(I_2+\sigma_2)} \\
q^{2(\sigma_2'-\sigma_2+1)},q^{2(2j(\lambda)+\sigma_2'-I_2-c+1)}
\end{array} q^2, q^{2(I_2+2j(\lambda)+\sigma_2'+2)}\right)\\
& \times q^{  
  j(\lambda)(2  I_2 -  2 \sigma_2') +c(- I_2   - 3    \sigma_2/2 +   \sigma_2'/2 ) - \sigma_2  \sigma_2'/2  + 3  \sigma_2^2/4 - 3  \sigma_2'^2/4+ I_2  \sigma_2'},
    \end{aligned}
\end{equation*}
and
\begin{equation*}  C_1(I_1,I_2,\sigma_2,\sigma_2')=\frac{(-1)^{1/2\sigma_2+1/2\sigma_2'-I_2}}{(q-q^{-1})^{1/2\sigma_2+1/2\sigma_2'-I_2}}
 q^{
I_2/2 + (3  \sigma_2)/4 - (5  \sigma_2')/4  + 2  I_1  \sigma_2   - I_2^2/2  }.
\end{equation*}}
{  Note that $C_1(I_1,I_2,\sigma_2,\sigma_2')$ does not depend on $c,\sigma_1,\sigma_1',\lambda$ and gives a constant  in $B_{ml}(z,\lambda)$ thus can be discarded in duality function.}

{
Next, we convert back to particle-hole notation and also apply various identities (e.g.(1.5.2),(3.2.4),(3.2.5) in \cite{GR04}) to simplify the results. For example, up to a constant,
\normalsize
\begin{equation*}
    \begin{aligned}
         \lim_{J_0\xrightarrow[]{}\infty}&\langle \mu | B_{2J_0-c+N(\xi-\mu),2J_0-c}(\infty,\lambda) | \xi \rangle\\
   =\prod_{i=1}^N&  \mathcal{F}(I_i,j(\lambda)+N_{[1,i-1]}(\mu-I),c+N_{[1,i-1]}(\mu-\xi),\xi_i-I_i,\mu_i-I_i)\\
   =\prod_{i=1}^N&\frac{[2j(\lambda)-c+N_{[1,i-1]}(\mu+\xi-2I)]![\mu_i-\xi_i]!^{-1}}{\left(q^{4\left(j(\lambda)+1+N_{[1,i]}(\mu-I)\right)};q^2\right)_{2I_i-\mu_i}[2j(\lambda)-c+N_{[1,i-1]}(\mu+\xi-2I)+\mu_i-2I_i]!}\\
  &\times\sqrt{\frac{[\mu_i]![2I_i-\xi_i]![c+N_{[1,i]}(\mu-\xi)]!^{-1}}{[2I_i-\mu_i]![\xi_i]![c+N_{[1,i-1]}(\mu-\xi)]!}} q^{   2I_iN_{[1,i-1]}(\mu)+I_i(\mu_i-\xi_i)-\xi_iN_{[1,i-1]}(\mu)}\\
  &\times{ }_3 \varphi_2\left(\begin{array}{c}
q^{-2\left(c+N_{[1,i-1]}(\mu-\xi)\right)}, q^{-2(2I_i-\mu_i)},q^{-2\xi_i} \\
q^{2(\mu_i-\xi_i+1)},q^{2(2j(\lambda)-c+N_{[1,i-1]}(\mu+\xi-2I)+\mu_i-2I_i+1)}\end{array} q^2, q^{2\left(2\left(j(\lambda)+N_{[1,i-1]}(\mu-I)\right)+\mu_i+2\right)}\right)\\
=\prod_{i=1}^N&\frac{[2j(\lambda)-c+N_{[1,i-1]}(\mu+\xi-2I)]!}{\left(q^{4\left(j(\lambda)+1+N_{[1,i]}(\mu-I)\right)};q^2\right)_{2I_i-\mu_i}[2j(\lambda)-c+N_{[1,i-1]}(\mu+\xi-2I)+\mu_i-2I_i]!} \\
  &\times\frac{1}{\sqrt{[\mu_i]![2I_i-\xi_i]![2I_i-\mu_i]![\xi_i]!}}q^{  
  2I_iN_{[1,i-1]}(\mu)+I_i(\mu_i+\xi_i)-\xi_iN_{[1,i]}(\mu)} \\
  &\times{ }_3 \varphi_2\left(\begin{array}{c}
 q^{-2(2I_i-\mu_i)},q^{-2\xi_i} ,q^{2(2j(\lambda)+2N_{[1,i-1]}(\mu-I)+\mu_i-2I_i+1)},\\
q^{-4I_i},q^{2(2j(\lambda)-c+N_{[1,i-1]}(\mu+\xi-2I)+\mu_i-2I_i+1)}\end{array} q^2, q^{2}\right)\\
=\prod_{i=1}^N&\frac{\left(q^{2(2j(\lambda)-c+N_{[1,i-1]}(\mu+\xi-2I)+\mu_i-2I_i+1)};q^2\right)_{2I_i-\mu_i}}{\left(q^{4\left(j(\lambda)+1+N_{[1,i]}(\mu-I)\right)};q^2\right)_{2I_i-\mu_i}\sqrt{[\mu_i]![2I_i-\xi_i]![2I_i-\mu_i]![\xi_i]!}}\\
  &\times{ }_3 \varphi_2\left(\begin{array}{c}
 q^{-2(2I_i-\mu_i)},q^{-2\xi_i} ,q^{2(2j(\lambda)+2N_{[1,i-1]}(\mu-I)+\mu_i-2I_i+1)},\\
q^{-4I_i},q^{2(2j(\lambda)-c+N_{[1,i-1]}(\mu+\xi-2I)+\mu_i-2I_i+1)}\end{array} q^2, q^{2}\right)\\
& \times q^{  
  2I_iN_{[1,i-1]}(\mu-\xi)+I_i(\mu_i+\xi_i)-2\xi_iN_{[1,i]}(\mu)}.\\
    \end{aligned}
\end{equation*}
Q.E.D.
}

{\subsubsection{Proof of Proposition \ref{prop: G} and \ref{prop: PG}} \label{sec: GS}}
We first  perform the calculation for $G$ defined in  \eqref{eq: GSformula}.
Note that $G$ is not unique, we solve for one such that  \eqref{eq: PG} is satisfied. Next, we show that  with the $G$ given in  Proposition \ref{prop: G}, \eqref{eq: PG} is satisfied with the constant  $C(I_1,I_2,\xi_1+\xi_2)=q^{2(I_1+I_2)(I_1+I_2-\xi_1-\xi_2)}$.

{
Setting $l=2J_0$ and taking the limit $J_0\to\infty$ in \eqref{eq: GSformula}, using the first limit found in Proposition \ref{prop: limits}, we find that
\begin{align*}
     &\langle \xi_1,\xi_{2} \vert G(\lambda) \vert \xi_1,\xi_{2} \rangle=  \lim_{J_0\xrightarrow[]{}\infty}\langle \xi | B_{2J_0-\xi_1-\xi_2,2J_0}(\infty, \lambda) | 0,0\rangle\\    
    &=\frac{\left([2I_1-\xi_1]![2I_2-\xi_{2}]![\xi_{2}]![\xi_1]!\right)^{-1/2}}{\left(q^{4(j(\lambda)+1+\xi_1-I_1)};q^2\right)_{2I_1-\xi_1}\left(q^{4(j(\lambda)+1+\xi_1+\xi_{2}-I_1-I_2)};q^2\right)_{2I_2-\xi_{2}}}q^{(I_1+2I_2)\xi_1+I_2\xi_2}.
\end{align*}

}

{Next, by direct calculation, we have   
\begin{equation*}
\frac{\langle \xi_1,\xi_{2} \vert \Pi^{\otimes 2} G(-\mathbf{i}\infty)\Pi^{\otimes 2} \vert \xi_1,\xi_{2} \rangle}{\langle \xi_1,\xi_{2} \vert P G(-\mathbf{i}\infty)P \vert \xi_1,\xi_{2} \rangle}  =q^{2(I_1+I_2)(I_1+I_2-\xi_1-\xi_2)}.
   \end{equation*}
Thus $\langle \mu_1,\mu_2|C|\xi_1,\xi_2\rangle=\delta_{\mu_i=\xi_i}q^{2(I_1+I_2)(I_1+I_2-\xi_1-\xi_2)}$ in \eqref{eq: PG}. 
}

\subsubsection{Proof of Proposition \ref{prop: multiG}}\label{sec: multiG}

{ Similar to the two site case, $B_{ml}(\lambda)  \vert \Vec{0}\rangle$ is an eigenvector of $\mathcal{T}(z,{\lambda})$,  so let 
\begin{equation*}
    \begin{aligned}
        &\langle \mu|\mathcal{G}(\lambda)|\mu\rangle
        =\lim_{J_0\rightarrow \infty}\langle\mu|B_{2J_0-N(\mu),2J_0}(\lambda)  \vert \Vec{0}\rangle\\
      &  =C(N(\mu))\prod_{i=1}^N\frac{q^{  
  2I_iN_{[1,i-1]}(\mu)+I_i\mu_i}}{\sqrt{[\mu_i]![2I_i-\mu_i]!}\left(e^{-2\pi \mathbf{i}\lambda}q^{4N_{[1,i]}(\mu-I)+2};q^2\right)_{2I_i-\mu_i}}.
    \end{aligned}
\end{equation*}
Then conjugating $\mathcal{T}(z,{\lambda})$ by $\mathcal{G}$ makes it stochastic, and we will  show it is equal to $ \mathcal{T}^{\mathrm{stoch}}$ by induction.
}

{
When $N=2$, 
$$\langle \mu|\mathcal{G}(\lambda)|\mu\rangle=\langle \mu|G_{1,2}(\lambda)|\mu\rangle,
$$
it is immediate that conjugating $\mathcal{T}(z,{\lambda})$ by $\mathcal{G}$ produces a stochastic matrix.
Now we assume there are $N$ sites, and define
\begin{equation*}
    \mathcal{T}_N^{\mathrm{stoch}}(z,{\lambda})=\mathcal{G}_N^{-1}(\lambda) \mathcal{T}_N(z,{\lambda})\mathcal{G}_N(\lambda) ,
\end{equation*}
 it suffices to show that 
\begin{multline*}
   \check{S}_{N,N+1} (z,\lambda-2\eta N_{[1,N-1]}(H))\mathcal{T}_N^{\mathrm{stoch}}(z,{\lambda})=\\
   \mathcal{G}_{N+1}^{-1}\check{R}_{N,N+1}(z,\lambda-2\eta N_{[1,N-1]}(H)) \mathcal{T}_N(z,{\lambda})\mathcal{G}_{N+1}.
\end{multline*}

}

{
Recalling that $\check{R}_{i,i+1}:V_1\otimes V_{i+1}\xrightarrow[]{}V_{i+1}\otimes V_1$, we obtain
\begin{equation*}
    \begin{aligned}
       & \frac{\langle \xi|G_{N,N+1}(\lambda-2\eta N_{[1,N-1]}(H))\check{R}_{N-1,N}|\mu\rangle}{\langle \xi|\check{R}_{N-1,N}G_{N,N+1}(\lambda-2\eta N_{[1,N-1]}(H))|\mu\rangle}\\
        &=\frac{\langle \xi_1,\xi_{N+1}|G(\lambda-2\eta N_{[1,N-1]}(H))|_{V_1\otimes V_{N+1}}|\xi_1,\xi_{N+1}\rangle}{\langle \mu_N,\mu_{N+1}|G(\lambda-2\eta N_{[1,N-1]}(H))|_{V_N\otimes V_{N+1}}|\mu_N,\mu_{N+1}\rangle}  \\
        &=\frac{\sqrt{[2I_N-\mu_N]![\mu_N]!}\left(e^{-2\pi \mathbf{i}\lambda}q^{4N_{[1,N]}(\mu-I)+2};q^2\right)_{2I_N-\mu_N}}{\sqrt{[2I_1-\xi_1]![\xi_1]!}\left(e^{-2\pi \mathbf{i}\lambda}q^{4N_{[1,N]}(\xi-J)+2};q^2\right)_{2I_1-\xi_1}}q^{(I_1+2I_{N+1})\xi_1-(I_N+2I_{N+1})\mu_N}.
    \end{aligned}
\end{equation*}
and
\begin{equation*}
    \begin{aligned}
       & \frac{\langle \xi|G_{N,N+1}(\lambda-2\eta N_{[1,N-1]}(H)) \mathcal{T}_N(z,{\lambda})\mathcal{G}_N(\lambda)|\mu\rangle}{\langle \xi| \mathcal{T}_N(z,{\lambda})\mathcal{G}_{N+1}(\lambda)|\mu\rangle}\\
        &=\frac{\langle \xi|G_{N,N+1}(\lambda-2\eta N_{[1,N-1]}(H)) \mathcal{T}_N(z,{\lambda})|\mu\rangle}{\langle \xi| \mathcal{T}_N(z,{\lambda})G_{N,N+1}(\lambda-2\eta N_{[1,N-1]}(H))|\mu\rangle}  \frac{\langle\mu|G_{N,N+1}(\lambda-2\eta N_{[1,N-1]}(H))\mathcal{G}_{N}(\lambda)|\mu\rangle}{\langle\mu|\mathcal{G}_{N+1}(\lambda)|\mu\rangle}\\
        &=\frac{1}{\sqrt{[2I_1-\xi_1]![\xi_1]!}}q^{-2I_{N+1}N_{[1,N]}(\xi)+(I_1+2I_{N+1})\xi_1}.
    \end{aligned}
\end{equation*}
Similarly, 
\begin{equation*}
    \begin{aligned}
       \frac{\langle \xi|\check{R}_{N,N+1}\mathcal{G}_N^{-1}|\mu\rangle}{\langle \xi|\mathcal{G}_N^{-1}\check{R}_{N,N+1}|\mu\rangle}
        &=\frac{\sqrt{[2I_1-\mu_1]![\mu_1]!}\left(e^{-2\pi \mathbf{i}\lambda}q^{4N_{[1,N]}(\mu-I)+2};q^2\right)_{2I_1-\mu_1}}{\sqrt{[2I_{N+1}-\xi_{N+1}]![\xi_{N+1}]!}\left(e^{-2\pi \mathbf{i}\lambda}q^{4N_{[2,N+1]}(\xi-J)+2};q^2\right)_{2I_{N+1}-\xi_{N+1}}}\\
        &q^{2(I_{N+1}-I_1)N_{[2,N]}(\xi)+I_{N+1}\xi_{N+1}-I_1\mu_1}.
    \end{aligned}
\end{equation*}
Thus,
\begin{equation*}
    \begin{aligned}
      &  \frac{\langle \xi|G^{-1}_{N,N+1}(\lambda-2\eta N_{[1,N-1]}(H)))\check{R}_{N,N+1}\mathcal{G}_{N}^{-1}|\mu\rangle}{\langle \xi|\mathcal{G}_{N+1}^{-1}\check{R}_{N,N+1}|\mu\rangle}\\
       &=\frac{\langle \xi|\check{R}_{N,N+1}\mathcal{G}_N^{-1}|\mu\rangle}{\langle \xi|\mathcal{G}_N^{-1}\check{R}_{N,N+1}|\mu\rangle}\frac{\langle \xi|G^{-1}_{N,N+1}(\lambda-2\eta N_{[1,N-1]}(H))\mathcal{G}^{-1}_N|\xi\rangle}{\langle \xi|\mathcal{G}^{-1}_{N+1}|\xi\rangle}\\
        &=\sqrt{[2I_1-\mu_1]![\mu_1]!}
        q^{2I_{N+1}N_{[2,N]}(\mu)-I_1\mu_1}.
    \end{aligned}
   \end{equation*}
Combining the above equations, we obtain the desired equality. Q.E.D.}

\subsubsection{Proof of Corollary \ref{cor: ort}}
{
This section provides explicit calculations for Corollary \ref{cor: ort}. 
}

{
Recall from \cite{GR04} that 
\begin{equation*}
{ }_3 \varphi_2\left[\begin{array}{c}
q^{-n}, a, b \\
d, e
\end{array} ; q, \frac{d e q^n}{a b}\right]=\frac{(e / a ; q)_n}{(e ; q)_n}{ }_3 \phi_2\left[\begin{array}{l}
q^{-n}, a, d / b \\
d, a q^{1-n} / e
\end{array} ; q, q\right]
\end{equation*}
and  the inversion formula
$$
\begin{aligned}
& { }_r \phi_s\left[\begin{array}{l}
a_1, \ldots, a_r \\
b_1, \ldots, b_s
\end{array} ;q, z\right] \\
& =\sum_{n=0}^{\infty} \frac{\left(a_1^{-1}, \ldots, a_r^{-1} ; q^{-1}\right)_n}{\left(q^{-1}, b_1^{-1}, \ldots, b_s^{-1} ; q^{-1}\right)_n}\left(\frac{a_1 \cdots a_r z}{b_1 \cdots b_s q}\right)^n .
\end{aligned}
$$
}

{
Combining the above two identities we get 
\begin{equation*}
  { }_3 \phi_2\left[\begin{array}{l}
q^{n}, 1/a, 1 / b \\
1/d, 1 / e
\end{array} ; q^{-1}, q^{-1}\right]  =\frac{(e / a ; q)_n}{(e ; q)_n}{ }_3 \phi_2\left[\begin{array}{l}
q^{-n}, a, d / b \\
d, a q^{1-n} / e
\end{array} ; q, q\right].
\end{equation*}
 Taking $e$ to $\infty$ yields
\begin{equation}\label{eq: trans}
  { }_3 \phi_2\left[\begin{array}{l}
q^{n}, 1/a, 1 / b \\
1/d, 0
\end{array} ; q^{-1}, q^{-1}\right]  =a^{-n}{ }_3 \phi_2\left[\begin{array}{l}
q^{-n}, a, d / b \\
d, 0
\end{array} ; q, q\right].
\end{equation}
}

{
 Applying \eqref{eq: trans} by taking $a=q^{-2(2I_i-\mu_i)}$, $n=\xi_i$, $d=q^{-4I_i}$ and $b^{-1}=e^{-2\pi \mathbf{i}\lambda}q^{2(2N_{[1,i-1]}(\mu-I)+\mu_i)}$ to the following result,
\begin{equation*}
    \begin{aligned}
  \lim_{c\xrightarrow[]{}-\infty} \mathcal{D}_{c}(\mu,\xi)
 =\prod_{i=1}^N &{ }_3 \varphi_2\left(\begin{array}{c}
 q^{-2(2I_i-\mu_i)},q^{-2\xi_i} ,e^{-2\pi \mathbf{i}\lambda}q^{2(2N_{[1,i-1]}(\mu-I)+\mu_i-2I_i)}\\
q^{-4I_i},0\end{array} q^2, q^{2}\right) \\
& \times q^{  -4I_iN_{[1,i-1]}(\xi)-2\xi_iN_{[1,i]}(\mu)}.\\
 \end{aligned}
\end{equation*}
 we obtain that 
\begin{equation*}
    \begin{aligned}
  \lim_{c\xrightarrow[]{}-\infty} \mathcal{D}_{c}(\mu,\xi)=\prod_{i=1}^N&{ }_3 \varphi_2\left(\begin{array}{c}
 q^{2(2I_i-\mu_i)},q^{2\xi_i} ,e^{-2\pi \mathbf{i}\lambda}q^{2(2N_{[1,i-1]}(\mu-I)+\mu_i)}\\
 q^{4I_i},0\end{array} q^{-2}, q^{-2}\right) \\&  \times q^{ -2\xi_i(2I_i-\mu_i) 
  -4I_iN_{[1,i-1]}(\xi)-2\xi_iN_{[1,i]}(\mu)}\\
  =\prod_{i=1}^N&{ }_3 \varphi_2\left(\begin{array}{c}
 q^{2(2I_i-\mu_i)},q^{2\xi_i} ,e^{-2\pi \mathbf{i}\lambda}q^{2(2N_{[1,i-1]}(\mu-I)+\mu_i)}\\
q^{4I_i},0\end{array} q^{-2}, q^{-2}\right) \\&  \times q^{ 
  -4I_iN_{[1,i]}(\xi)-2\xi_iN_{[1,i-1]}(\mu)},
 \end{aligned}
\end{equation*}
which is equal to $\mathcal{D}_{ort}(\mu,\xi)$ up to a constant. Q.E.D.
}

\subsection{Proof of Proposition \ref{prop: matching}}
{
In \cite{AA18}, the author  computes explicit formulas for the stochastic weights of the dynamic stochastic higher spin vertex model by using fusion. Essentially, he proves a recursive relation for the weights and then solves the relation explicitly. Here, we prove that the same recursive relation is satisfied by our weights. Because the (difficult) task of solving the weights was already done, a proof of the recursive relation suffices to prove Proposition \ref{prop: matching}. While ``spiritually'' the fusion of Aggarwal is the same as the algebraic fusion that arises from the universal $R$-matrix and twister, for reasons of completeness a proof still needs to be provided.
}

{
The first step is to translate between the notational conventions used by different authors. Since this is most naturally illustrated through examples, we defer the discussion to the Appendix. In this section, we prove that the stochastic matrix $S$ has entries given by \eqref{eq: Amolweights}.

}

{
Recall from section \ref{QGas} that for each positive integer,  $V^{2J}$ denotes the $(2J+1)$-dimensional representation of \(U_q( \mathfrak{sl}_2)\). There is a projection 
$$
V^{1} \otimes V^{J-1} \rightarrow \mathcal{P}\left( V^{1} \otimes V^{J-1}\right) \cong V^J,
$$
where $\mathcal{P}$ is the projection $\frac{1}{q-q^{-1}}\check{R}^{J-1,1}(q^2,\lambda)$. One can verify by direct calculation that $\mathcal{P}^2 = \mathcal{P}$ using the explicit formulas for the $R$-matrix.  Let $\iota:V^J\rightarrow V^{J-1} \otimes V^1$ denote the inclusion.
}

{
Below, we use the same notations as \cite{AA18}. By the dynamic Yang-Baxter Equation \eqref{YBE}, where we let $R_{13}$ be $W_{I}$ and $R_{23}$ be $\widehat{W}_{J-1}$, we see that
\begin{multline*}
\langle j_2,i_2 \vert W_J \vert j_1,i_1 \rangle = \sum_{k=0}^1  \sum_{l=0}^1 \sum_{m=0}^1 \langle j_2,i_2 \vert \mathcal{P} \otimes \mathrm{id} \vert m,l,i_2\rangle \langle m,l, i_2 \vert W_I \vert k,l,\delta\rangle\\
\times \langle k,l  ,  \delta \vert \widehat{W}_{J-1} \vert k,j_1-k, i_1 \rangle \langle k,j_1-k,i_1\vert \iota \otimes \mathrm{id} \vert j_1,i_1 \rangle  
\end{multline*}
where $l+\delta = j_1-k+i_1$. Note that this term is only nonzero if $j_2=m+l$, so the summation only has four terms. In the above equation, the dependence on the parameters is dropped. 
}

{Having set $u=q^{2(J-1)}$ in \eqref{YBE} so that $\mathcal{P}=R_{12}$, the two terms are now 
$$
W_1(i_2,m;\delta,k\vert q^{2(J-1)}z,\lambda), \quad \quad \widehat{W}_{J-1}( i_1,j_1-k;\delta,l\vert z,\lambda -2\eta H_1(k)).
$$
We write $\delta$ in terms of $i_2,m,k$ by 
$$
\delta = j_1+i_1-k-l = j_2+i_2 - k - l = i_2-k+m.
$$
We can see that this matches with the four terms in the recurrence relation in \cite{AA18};
$$
\begin{array}{l} \widehat{W}_{J}\left(i_{1}, j_{1} ;\right.\left.i_{2}, j_{2} \mid v, \lambda\right) \\ =\widehat{W}_{J-1}\left(i_{1}, j_{1} ; i_{2}, j_{2} \mid v, \lambda-2 \eta\right) W_{1}\left(i_{2}, 0 ; i_{2}, 0 \mid v+2 \eta(J-1), \lambda\right) \\ +\widehat{W}_{J-1}\left(i_{1}, j_{1}-1 ; i_{2}-1, j_{2} \mid v, \lambda+2 \eta\right) W_{1}\left(i_{2}-1,1 ; i_{2}, 0 \mid v+2 \eta(J-1), \lambda\right) \\ +\widehat{W}_{J-1}\left(i_{1}, j_{1} ; i_{2}+1, j_{2}-1 \mid v, \lambda-2 \eta\right) W_{1}\left(i_{2}+1,0 ; i_{2}, 1 \mid v+2 \eta(J-1), \lambda\right) \\ +\widehat{W}_{J-1}\left(i_{1}, j_{1}-1 ; i_{2}, j_{2}-1 \mid v, \lambda+2 \eta\right) W_{1}\left(i_{2}, 1 ; i_{2}, 1 \mid v+2 \eta(J-1), \lambda\right)\end{array}
$$
where $e^{v}=z$. Therefore the weights match and Proposition \ref{prop: matching} is proved. 
}

\section{Appendix: Matching of weights }
In this appendix, we briefly explain how to identify the weights \eqref{eq: Amolweights} with the stochastic matrix $S(z,\lambda)$. For simplicity, we restrict attention to the case of spectral parameter $u=0$ in \eqref{eq: Amolweights} and match only the hypergeometric terms.

{ First, note that
\normalsize
\begin{equation*}
\begin{aligned}
&\lim_{u\xrightarrow[]{}0} \psi\left(i_1, j_1 ; i_2, j_2\right) \\
& =\mathbf{1}_{i_1+j_1=i_2+j_2} \left(\frac{q^{j_1 - 2i_1j_1 - 4i_2j_1 - 2i_1^2 + 2i_2^2 + j_1^2}}{s^{2j_1+2j_2}\kappa^{j_2}}\right) \frac{\left(q^{2i_1-2j_2+2} ; q^2\right)_{j_2}\left(q^{2j_2-2J} ; q^2\right)_{j_1}\left(s^2 q^{2i_1-2j_2} ; q^2\right)_{j_1}}{(q^2 ; q^2)_{j_2}\left(q^{2j_2-2J} ; q^2\right)_{j_1-j_2}} \\
&  \times  \frac{\left(q^2 \kappa^{-1} ; q^2\right)_{j_1}\left(q^{2j_2-4 i_1+2} s^{-2} \kappa^{-1} ; q^2\right)_{i_1-j_2}}{\left(q^{2-2j_2} \kappa^{-1} ; q^2\right)_{j_1}\left(q^{2j_2-2i_1-2J+2} s^{-2} \kappa^{-1} ; q^2\right)_{j_1}\left(q^{2j_2-4 i_1-2J} s^{-2} \kappa^{-1} ; q^2\right)_{j_2}\left(q^{4 j_2-4 i_1-2J+2} s^{-2} \kappa^{-1} ; q^2\right)_{i_1-j_2}} \\
&  \times{ }_{8} W_7\left(q^{-2j_2} \kappa^{-1} ; q^{-2j_1}, q^{-2j_2}, q^{2j_1-2J} \kappa^{-1}, q^{2-2i_1} s^{-2} \kappa^{-1}, q^{-2i_1} \kappa^{-1};q^2,s^2\kappa q^{2+4i_1-2j_2+2J}\right) 
\end{aligned}
\end{equation*}
}
{

Now  set $q^{2j(\lambda)+1}=e^{-2\pi\mathbf{i}\lambda}$ and     apply   \cite[(2.5.1) and (3.2.1)]{GR04} to get
\normalsize
\begin{equation*}
    \begin{aligned}
        &{ }_{8} W_7\left(q^{-2j_2} \kappa^{-1} ; q^{-2j_1}, q^{-2j_2}, q^{2j_1-2J} \kappa^{-1}, q^{2-2i_1} s^{-2} \kappa^{-1}, q^{-2i_1} \kappa^{-1};q^2,s^2\kappa q^{2+4i_1-2j_2+2J}\right)  \\
        =&\frac{\left(q^{2J - 2i_1 + 2i_2 - 6  j_1 - 4  j(\lambda)},q^{2i_2 + 4  j(\lambda) + 4}\right)_{i_1 -  I- J + 2  j_1 + 2  j(\lambda)}\left(q^{-2j_1},q^{2J +  2\Lambda - 2i_2 - 2j_1 -4 j(\lambda) + 2}\right)_{J +  I - i_1 - 2  j_1 - 2  j(\lambda)}}{\left(q^{2J -2 i_1 +2 i_2 - 4  j_1 + 2},q^{ 2i_2 - 2j_1 + 1}\right)_{i_1 -  I - J + 2  j_1 + 2  j(\lambda)}\left(q^{2J -4  j_1 - 4  j(\lambda)},q^{ 2I - 2i_2 + 2}\right)_{J +  I - i_1 - 2  j_1 - 2  j(\lambda)}}\\
        &\times { }_4 \varphi_3\left(\begin{array}{c}
q^{2(J +  I - i_1 - 2  j_1 - 2  j(\lambda))}, q^{2(J - j_1 - 2  j(\lambda))}, q^{2(J +    I - i_1 - i_2 - j_1 - 2  j(\lambda))},q^{2(I - i_2 - 2  j(\lambda))} \\
q^{2(J +    I - i_1 - i_2 - 2  j_1 - 4  j(\lambda) - 1)},q^{2(J +    I - i_2 - j_1 - 2  j(\lambda) + 1)}, q^{2( J +    I - i_1 - j_1 - 2  j(\lambda) + 1)} 
\end{array} q^{2}, q^{2}\right)
    \end{aligned}
\end{equation*}
}

{
On the other hand, applying a symmetry of $6j$--symbol (see \cite{KR88}) and the equivalent definition through  Racah--Wilson polynomial \eqref{eq: 6j vs RW}, we have 
\normalsize
\begin{equation*}
\begin{aligned}
\normalsize
 &\langle i_1,j_1| S^{I,J}(\infty,\lambda)| i_2, j_2\rangle=\langle j_1,i_1 \vert G(\lambda)^{-1} \vert j_1,i_1 \rangle\langle i_1,j_1| R^{I,J}(\infty,\lambda)| i_2, j_2\rangle \langle i_2, j_2 \vert G(\lambda) \vert i_2, j_2\rangle\\&\sim\left\{\begin{array}{ccc}
I/2 &j(\lambda)+j_1+i_1-J/2-  I/2 &  j(\lambda)+j_1-J/2\\
J/2& j(\lambda) & j(\lambda)+i_2-  \Lambda/2
\end{array}\right\}_q\\
&\sim \left\{\begin{array}{ccc}
i_2/2 -   I/2 + j_2/2 + j(\lambda) &i_1 -  I/2 - J/2 + j_2 + j(\lambda)& J/2 + i_2/2 - j_2/2 \\
i_2/2 - J/2 + j_2/2 + j(\lambda)& j(\lambda) & I/2 - i_2/2 + j_2/2
\end{array}\right\}_q\\
&\sim { }_4 \varphi_3\left(\begin{array}{c}
q^{2(J +  I - i_1 - 2  j_1 - 2  j(\lambda))}, q^{2(J - j_1 - 2  j(\lambda))}, q^{2(J +    I - i_1 - i_2 - j_1 - 2  j(\lambda))},q^{2(I - i_2 - 2  j(\lambda))} \\
q^{2(J +   I - i_1 - i_2 - 2  j_1 - 4  j(\lambda) - 1)},q^{2(J +    I - i_2 - j_1 - 2  j(\lambda) + 1)}, q^{2( J +    I - i_1 - j_1 - 2  j(\lambda) + 1)} 
\end{array} q^{2}, q^{2}\right)\\
&\sim\lim_{u\xrightarrow[]{}0} \psi\left(i_1, j_1 ; i_2, j_2\right)
\end{aligned}
\end{equation*}

where $A\sim B$ means that the hypergeometric term ${}_4\varphi_3$ appearing in $A$ coincides with that appearing in $B$, although the prefactors may differ.

}

\bibliographystyle{plain}
\bibliography{mainMay2024}

@article{AA18,
  author    = {Amol Aggarwal},
  title     = {Dynamical stochastic higher spin vertex models},
  journal   = {Selecta Math. (N.S.)},
  volume    = {24},
  number    = {3},
  pages     = {2659--2735},
  year      = {2018},
}

@article{AW79,
  author    = {Richard Askey and James Wilson},
  title     = {A set of orthogonal polynomials that generalize the {R}acah coefficients or {$6-j$} symbols},
  journal   = {SIAM J. Math. Anal.},
  volume    = {10},
  number    = {5},
  pages     = {1008--1016},
  year      = {1979},
}

@article{ACR18,
  author    = {Mario Ayala and Gioia Carinci and Frank Redig},
  title     = {Quantitative {B}oltzmann-{G}ibbs principles via orthogonal polynomial duality},
  journal   = {J. Stat. Phys.},
  volume    = {171},
  number    = {6},
  pages     = {980--999},
  year      = {2018},
}

@article{ACR21,
  author    = {Mario Ayala and Gioia Carinci and Frank Redig},
  title     = {Higher order fluctuation fields and orthogonal duality polynomials},
  journal   = {Electron. J. Probab.},
  volume    = {26},
  pages     = {Paper No. 27, 35},
  year      = {2021},
}

@article{BBB95,
  author    = {O. Babelon and D. Bernard and E. Billey},
  title     = {A quasi-{H}opf algebra interpretation of quantum {$3$}-{$j$} and {$6$}-{$j$} symbols and difference equations},
  journal   = {Phys. Lett. B},
  volume    = {375},
  number    = {1-4},
  pages     = {89--97},
  year      = {1996},
}

@article{BOR17,
  author    = {Alexei Borodin},
  title     = {Symmetric elliptic functions, {IRF} models, and dynamic exclusion processes},
  journal   = {J. Eur. Math. Soc. (JEMS)},
  volume    = {22},
  number    = {5},
  pages     = {1353--1421},
  year      = {2020},
}

@article{BC20,
  author    = {Alexei Borodin and Ivan Corwin},
  title     = {Dynamic {ASEP}, duality, and continuous {$q^{-1}$}-{H}ermite polynomials},
  journal   = {Int. Math. Res. Not. IMRN},
  number    = {3},
  pages     = {641--668},
  year      = {2020},
}

@article{BCG16,
  author    = {Alexei Borodin and Ivan Corwin and Vadim Gorin},
  title     = {Stochastic six-vertex model},
  journal   = {Duke Math. J.},
  volume    = {165},
  number    = {3},
  pages     = {563--624},
  year      = {2016},
}

@article{CFGR18,
  author    = {Gioia Carinci and Chiara Franceschini and Cristian Giardin\`{a} and Wolter Groenevelt and Frank Redig},
  title     = {Orthogonal dualities of {M}arkov processes and unitary symmetries},
  journal   = {SIGMA Symmetry Integrability Geom. Methods Appl.},
  volume    = {15},
  pages     = {Paper No. 053, 27},
  year      = {2019},
}

@article{CFG21,
  author    = {Gioia Carinci and Chiara Franceschini and Wolter Groenevelt},
  title     = {{$q$}-{O}rthogonal dualities for asymmetric particle systems},
  journal   = {Electron. J. Probab.},
  volume    = {26},
  pages     = {Paper No. 108, 38},
  year      = {2021},
}

@article{CGR94,
  author    = {Eug\`{e}ne Cremmer and Jean-Loup Gervais and Jean-Fran\c{c}ois Roussel},
  title     = {The quantum group structure of {$2$}{D} gravity and minimal models. {II}. {T}he genus-zero chiral bootstrap},
  journal   = {Comm. Math. Phys.},
  volume    = {161},
  number    = {3},
  pages     = {597--630},
  year      = {1994},
}

@article{FV96,
  author    = {Giovanni Felder and Alexander Varchenko},
  title     = {On representations of the elliptic quantum group {$E_{\tau,\eta}({\rm sl}_2)$}},
  journal   = {Comm. Math. Phys.},
  volume    = {181},
  number    = {3},
  pages     = {741--761},
  year      = {1996},
}

@misc{FKZ22,
  author    = {Chiara Franceschini and Jeffrey Kuan and Zhengye Zhou},
  title     = {Orthogonal polynomial duality and unitary symmetries of multi--species {ASEP}$(q,\boldsymbol{\theta})$ and higher--spin vertex models via $^*$--bialgebra structure of higher rank quantum groups},
  year      = {2022},
}

@book{GR04,
  author    = {George Gasper and Mizan Rahman},
  title     = {Basic hypergeometric series},
  series    = {Encyclopedia of Mathematics and its Applications},
  volume    = {96},
  edition   = {second},
  publisher = {Cambridge University Press},
  address   = {Cambridge},
  year      = {2004},
  note      = {With a foreword by Richard Askey},
}

@article{Wolter19,
  author    = {Wolter Groenevelt},
  title     = {Orthogonal stochastic duality functions from {L}ie algebra representations},
  journal   = {J. Stat. Phys.},
  volume    = {174},
  number    = {1},
  pages     = {97--119},
  year      = {2019},
}

@article{groenevelt2023generalized,
  author  = {Groenevelt, Wolter and Wagenaar, Carel},
     title = {A generalized dynamic asymmetric exclusion process: orthogonal
              dualities and degenerations},
   journal = {J. Phys. A},
  fjournal = {Journal of Physics. A. Mathematical and Theoretical},
    volume = {57},
      YEAR = {2024},
    NUMBER = {37},
     PAGES = {Paper No. 375202, 66},
      ISSN = {1751-8113,1751-8121},
   MRCLASS = {60K35 (81R50)},
  MRNUMBER = {4796968},
       DOI = {10.1088/1751-8121/ad6f7b},
       URL = {https://doi.org/10.1088/1751-8121/ad6f7b},
}

@article{JK14,
  author    = {Sabine Jansen and Noemi Kurt},
  title     = {On the notion(s) of duality for {M}arkov processes},
  journal   = {Probab. Surv.},
  volume    = {11},
  pages     = {59--120},
  year      = {2014},
}

@article{JKOS97,
  author    = {M. Jimbo and H. Konno and S. Odake and J. Shiraishi},
  title     = {Quasi-{H}opf twistors for elliptic quantum groups},
  journal   = {Transform. Groups},
  volume    = {4},
  number    = {4},
  pages     = {303--327},
  year      = {1999},
}

@incollection{KR88,
  author    = {A. N. Kirillov and N. Yu. Reshetikhin},
  title     = {Representations of the algebra {${U}_q({\rm sl}(2)),\;q$}-orthogonal polynomials and invariants of links},
  booktitle = {Infinite-dimensional {L}ie algebras and groups ({L}uminy-{M}arseille, 1988)},
  series    = {Adv. Ser. Math. Phys.},
  volume    = {7},
  pages     = {285--339},
  publisher = {World Sci. Publ.},
  address   = {Teaneck, NJ},
  year      = {1989},
}

@book{KSBOOK,
  author    = {Leonid I. Korogodski and Yan S. Soibelman},
  title     = {Algebras of functions on quantum groups. {P}art {I}},
  series    = {Mathematical Surveys and Monographs},
  volume    = {56},
  publisher = {American Mathematical Society},
  address   = {Providence, RI},
  year      = {1998},
}

@article{Kuan17,
  author    = {Jeffrey Kuan},
  title     = {An algebraic construction of duality functions for the stochastic {${U}_q(A_n^{(1)})$} vertex model and its degenerations},
  journal   = {Comm. Math. Phys.},
  volume    = {359},
  number    = {1},
  pages     = {121--187},
  year      = {2018},
}

@article{Kuan18,
  author    = {Jeffrey Kuan},
  title     = {A multi-species {${\rm ASEP}(q,j)$} and {$q$}-{TAZRP} with stochastic duality},
  journal   = {Int. Math. Res. Not. IMRN},
  number    = {17},
  pages     = {5378--5416},
  year      = {2018},
}

@misc{KZ23,
  author    = {Jeffrey Kuan and Zhengye Zhou},
  title     = {Asymptotics of two-point correlations in the multi-species q-{TAZRP}},
  year      = {2023},
}

@article{KMMO16,
  author    = {A. Kuniba and V. V. Mangazeev and S. Maruyama and M. Okado},
  title     = {Stochastic {$R$} matrix for {$U_q(A_n^{(1)})$}},
  journal   = {Nuclear Phys. B},
  volume    = {913},
  pages     = {248--277},
  year      = {2016},
}

@article{SS94,
  author    = {Gunter Sch\"{u}tz and Sven Sandow},
  title     = {Non-abelian symmetries of stochastic processes: Derivation of correlation functions for random-vertex models and disordered-interacting-particle systems},
  journal   = {Phys. Rev. E},
  volume    = {49},
  pages     = {2726--2741},
  month     = {Apr},
  year      = {1994},
}

@article{Schutz97,
  author    = {Gunter M. Sch\"{u}tz},
  title     = {Duality relations for asymmetric exclusion processes},
  journal   = {J. Statist. Phys.},
  volume    = {86},
  number    = {5-6},
  pages     = {1265--1287},
  year      = {1997},
}

@article{Zhou_2021,
  author    = {Zhengye Zhou},
  title     = {Orthogonal polynomial stochastic duality functions for multi-species {SEP}(2j) and multi-species {IRW}},
  journal   = {Symmetry, Integrability and Geometry: Methods and Applications},
  month     = {dec},
  year      = {2021},
}

}
\end{document}